# Universal pointwise selection rule in multivariate function estimation

ALEXANDER GOLDENSHLUGER[1] and OLEG LEPSKI[2]

[1]*Department of Statistics, University of Haifa, 31905 Haifa, Israel. E-mail: goldensh@stat.haifa.ac.il*
[2]*Laboratoire d'Analyse, Topologie et Probabilités UMR CNRS 6632, Université de Provence, 39, rue F. Joliot Curie, 13453 Marseille, France. E-mail: lepski@cmi.univ-mrs.fr*

In this paper, we study the problem of pointwise estimation of a multivariate function. We develop a general pointwise estimation procedure that is based on selection of estimators from a large parameterized collection. An upper bound on the pointwise risk is established and it is shown that the proposed selection procedure specialized for different collections of estimators leads to minimax and adaptive minimax estimators in various settings.

## 1. Introduction

In this paper, we study the problem of pointwise nonparametric estimation of an unknown function $F: \mathbb{R}^d \to \mathbb{R}$ in the multidimensional Gaussian white noise model

$$Y(\mathrm{d}t) = F(t)\,\mathrm{d}t + \varepsilon W(\mathrm{d}t), \qquad t = (t_1, \ldots, t_d) \in \mathcal{D}, \tag{1}$$

where $\mathcal{D}$ is an open interval in $\mathbb{R}^d$ containing $\mathcal{D}_0 := [-1/2, 1/2]^d$, $W$ is the standard Wiener process in $\mathbb{R}^d$ and $0 < \varepsilon < 1$ is the noise level. Our goal is to estimate $F$ at a given point $x \in \mathcal{D}_0$ using the observation $\mathcal{Y}_\varepsilon := \{Y(t), t \in \mathcal{D}\}$. We assume that the observation set $\mathcal{D}$ is larger than $\mathcal{D}_0$ in order to avoid boundary effects. Such assumptions are rather common in multivariate nonparametric models (see, e.g., Chen (1991), Hall (1989)).

Accuracy of an estimator $\tilde{F}(x) = \tilde{F}(x; \mathcal{Y}_\varepsilon)$ is measured by the risk

$$\mathcal{R}_r[\tilde{F}; F] := \{\mathbb{E}_F |\tilde{F}(x) - F(x)|^r\}^{1/r}, \qquad r > 0,$$

where $\mathbb{E}_F$ denotes the expectation with respect to the distribution $\mathbb{P}_F$ of $\mathcal{Y}_\varepsilon$ satisfying (1).

We develop a pointwise estimation procedure that is based on the selection of estimators from a large collection.

Denote by $\mathfrak{K}$ the set of all *kernels*, that is, functions $K: \mathcal{D} \times \mathcal{D}_0 \to \mathbb{R}$ such that $\int_\mathcal{D} K(t, x)\,\mathrm{d}t = 1$ for all $x \in \mathcal{D}_0$. Let $\mathcal{K}$ be a given subset of $\mathfrak{K}$ and let $\mathcal{F}(\mathcal{K})$ be the corresponding collection of linear estimators of $F(x)$ associated with the family $\mathcal{K}$:

$$\mathcal{F}(\mathcal{K}) := \left\{ \hat{F}^K(x) = \int_\mathcal{D} K(t,x) Y(\mathrm{d}t), K \in \mathcal{K} \right\}. \tag{2}$$

In this paper, we propose an estimator of $F(x)$ that is based on *random* (measurable with respect to $\mathcal{Y}_\varepsilon$) selection from the collection $\mathcal{F}(\mathcal{K})$. Denoting this estimator by $\hat{F}_\mathcal{K}(x)$, we have

$$\hat{F}_\mathcal{K}(x) = \hat{F}^{\hat{K}}(x),$$

where $\hat{K} \in \mathcal{K}$ for any "frozen" trajectory $\mathcal{Y}_\varepsilon$. Although $\hat{F}_\mathcal{K}(x)$ can be constructed for *any* $\mathcal{K} \in \mathfrak{K}$, we establish the upper bound on its risk only for $\mathcal{K} \in \mathcal{P}(\mathfrak{K})$; here, $\mathcal{P}(\mathfrak{K})$ is the set of all collections $\mathcal{K}$ satisfying some natural and non-restrictive conditions (see (K0)–(K2) in Section 2). We then prove (Theorem 1) that for all $\varepsilon$ small enough and for any $\mathcal{K} \in \mathcal{P}(\mathfrak{K})$,

$$\mathcal{R}_r[\hat{F}_\mathcal{K}; F] \leq \mathcal{U}_{\mathcal{K},F}(x) \qquad \forall F \in \mathbb{F}(\mathcal{K}), \tag{3}$$

where the upper bound $\mathcal{U}_{\mathcal{K},F}(x)$ is completely determined by the function $F$ and by the family of kernels $\mathcal{K}$. Here, $\mathbb{F}(\mathcal{K})$ is a large nonparametric set whose dependence on $\mathcal{K}$ is typically weak. In particular, in most interesting examples, we have $\mathbb{F}(\mathcal{K}) \supset \mathbb{C}_b(\mathcal{D})$ (see Remark 6 and Theorem 1 below), where $\mathbb{C}_b(\mathcal{D})$ is the set of all uniformly bounded continuous functions.

It is important to emphasize that our selection procedure can be applied to different collections of kernel estimators. Thus, we derive estimators $\{\hat{F}_\mathcal{K}, \mathcal{K} \in \mathcal{P}(\mathfrak{K})\}$ with different statistical properties as an output of a unique computational routine.

*Kernel collections.* Consider several examples of kernel collections for which the upper bound (3) can be established. Here and later on, $K: \mathbb{R}^d \to \mathbb{R}$ is a fixed function and for all $u, v \in \mathbb{R}^d$, we understand $u/v$ as $(u_1/v_1, \ldots, u_d/v_d)$.

**Example 1.** Let $d = 1$ and for any $x \in \mathcal{D}_0$, let

$$\mathcal{K}_1 = \left\{ h^{-1} K\left(\frac{\cdot - x}{h}\right), h \in [h_{\min}, h_{\max}] \right\},$$

where $0 < h_{\min} < h_{\max} \leq 1$ are given real numbers.

A random choice from this collection leading to a data-driven bandwidth $\hat{h}_\varepsilon(x) = \hat{h}(x, \mathcal{Y}_\varepsilon)$ was proposed in Lepski *et al.* (1997). The upper bound of type (3) obtained in that paper was used in order to establish minimax results on the Besov classes of functions. We note that the estimator $\hat{F}^{\hat{K}}, \hat{K}(t,x) = \hat{h}_\varepsilon(x)^{-1} K([t-x]/\hat{h}_\varepsilon(x))$ constructed in Lepski *et al.* (1997) and the estimator $\hat{F}_{\mathcal{K}_1}$ developed in this paper are different.

***Example 2.*** Consider generalization of the above collection $\mathcal{K}_1$ to an arbitrary dimension $d>1$. Let $\mathcal{H} = \bigotimes_{i=1}^{d}[h_{\min}^{(i)}, h_{\max}^{(i)}]$ and

$$\mathcal{K}_{\mathcal{H}} = \left\{ \left[\prod_{i=1}^{d} h_i^{-1}\right] K\left(\frac{\cdot - x}{h}\right), h \in \mathcal{H} \right\}. \tag{4}$$

A sophisticated random choice from the collection $\mathcal{K}_{\mathcal{H}}$ was proposed in Kerkyacharian *et al.* (2001). The corresponding upper bound of the type (3) allowed minimax results to be obtained on the anisotropic Besov classes of functions (functions with inhomogeneous smoothness). Again, we note that the estimator constructed in Kerkyacharian *et al.* (2001) and the estimator $\hat{F}_{\mathcal{K}_{\mathcal{H}}}$ proposed in the present paper are different.

Even though $\mathcal{K}_{\mathcal{H}}$ is a rather rich collection, it is not "sufficiently" rich for many interesting statistical problems. The next example illustrates this point.

***Example 3.*** Denote by $\mathcal{E}$ the set of all $d \times d$ orthogonal matrices and let

$$\mathcal{H}_1 = \{h \in \mathbb{R}^d : h = (h_1, h_{\max}, \ldots, h_{\max}), h_1 \in [h_{\min}, h_{\max}]\}.$$

Consider the kernel collection

$$\mathcal{K}_{SI} = \left\{ \frac{1}{h_1} \frac{1}{h_{\max}^{d-1}} K\left(\frac{E^T[\cdot - x]}{h}\right), E \in \mathcal{E}, h \in \mathcal{H}_1 \right\}. \tag{5}$$

This collection is appropriate for the estimation of functions possessing the single index structure. We refer to Chen (1991), Golubev (1992), Hristache *et al.* (2001) and references therein for works on estimation in the single index model.

Note that the collections (4) and (5) are quite different and "incomparable". However, one can easily define a more general collection of kernels that combines (4) and (5).

***Example 4.*** Define

$$K_{\mathcal{H},\mathcal{E}} = \left\{ \left[\prod_{i=1}^{d} h_i^{-1}\right] K\left(\frac{E^T(\cdot - x)}{h}\right), h \in \mathcal{H}, E \in \mathcal{E} \right\}.$$

The estimator $\hat{F}_{K_{\mathcal{H},\mathcal{E}}}$ could be applied simultaneously to estimate functions with inhomogeneous or unknown smoothness as well as functions with the single index structure.

The list of examples of kernel collections corresponding to different "structural" models (see Stone (1985)) could be continued. Selection from such collections leads to estimators that adapt simultaneously to a wide spectrum of assumptions on smoothness, structure, etc. Pointwise adaptive estimators based on selection from specific collections of estimators were also constructed in Lepski (1990, 1991), Lepski and Spokoiny (1997),

Goldenshluger and Nemirovski (1997), Tsybakov (1998), Klemelä and Tsybakov (2001) and Golubev (2004). A detailed discussion of relationships between our results and results in the cited papers is given in Section 3.3.

*Objective of the paper.* The local inequality (3) specialized to different families of kernels $\mathcal{K} \in \mathcal{P}(\mathfrak{K})$ allows us to derive minimax and adaptive results in various settings. This is the feature that characterizes the power of the estimator $\hat{F}_\mathcal{K}$ and usefulness of the upper bound in (3). In order to demonstrate universality of our selection procedure, we discuss its application to the following nonparametric estimation problems.

   (i) *Pointwise adaptive estimation in the single index model.* Here, we assume that $F(t) = f(\omega^T t)$, where $f : \mathbb{R} \to \mathbb{R}$ is an unknown function, $\omega \in \mathbb{S}^{d-1}$ is an unknown direction vector and $\mathbb{S}^{d-1}$ is the unit sphere in $\mathbb{R}^d$. Suppose, also, that $f$ belongs to the one-dimensional Hölder ball with unknown parameters. The objective is to estimate $F$ at a single given point $x \in \mathcal{D}_0$.
   (ii) *Pointwise minimax estimation over a union of anisotropic Hölder classes.* In this setting, it is assumed that $F$ belongs to the union of anisotropic Hölder classes $\mathbb{H}_d(\boldsymbol{\alpha}, L)$ (see Definition 2) over all $\boldsymbol{\alpha} = (\alpha_1, \ldots, \alpha_d)$ satisfying $\sum_{i=1}^d 1/\alpha_i = 1/\gamma$, where $\gamma > 0$ is a given number. The objective is to estimate $F(x)$ at a given point $x \in \mathcal{D}_0$.
   (iii) *Global minimax estimation over isotropic Besov classes.* Assume that $F$ belongs to the isotropic Besov class. The objective is to estimate $F$ globally on $\mathcal{D}_0$ with small $\mathbb{L}_r$-risk, $\mathcal{R}_{\mathbb{L}_r}[\tilde{F}; F]$, $r \in [1, \infty)$.

We are not aware of any results on problem (i) reported in the literature. For this problem, our procedure provides a minimax adaptive estimator in the sense of (6) with $\mathbb{F}_s^*$ being the one-dimensional Hölder class $\mathbb{H}_1(\alpha, L)$ and the parameter $s = (\alpha, L)$ including smoothness index $\alpha$ and constant $L$. Thus, in the setup of problem (i), there is no price to pay for adaptation to the unknown smoothness parameters $\alpha$ and $L$.

Problems (ii) and (iii) were considered in Klutchnikoff (2005) and Kerkyacharian *et al.* (2001), respectively. We note, however, that the methods proposed in these papers are highly specialized and are tailored to the problem in question. In contrast to this, we arrive at the solution to these problems by applying the same general selection procedure for different collections of estimators. In particular, our selection procedure applied to the collection $\mathcal{F}(\mathcal{K}_{SI})$ provides a solution to problem (i). The minimax estimators for problems (ii) and (iii) are constructed by using the proposed scheme on certain subcollections of $\mathcal{F}(\mathcal{K}_{\mathcal{H},\mathcal{E}})$. Moreover, we show that all of the problems (i)–(iii) can be solved simultaneously by the same selection procedure applied to the collection of estimators $\mathcal{F}(\mathcal{K}_{\mathcal{H},\mathcal{E}})$.

*Derivation of minimax and adaptive results.* Let us briefly discuss how to derive minimax and adaptive results from local inequalities of type (3).

In the framework of the minimax approach, $F$ is assumed to belong to some given set $\mathbb{F}^*$. The objective is to find an estimator $\hat{F}$ such that

$$\sup_{F \in \mathbb{F}^*} \mathcal{R}_r[\hat{F}; F] \asymp \inf_{\tilde{F}} \sup_{F \in \mathbb{F}^*} \mathcal{R}_r[\tilde{F}; F] \qquad \text{as } \varepsilon \to 0,$$

where inf is taken over all possible estimators. Here and in what follows, $a \asymp b$ means that $0 < c_1 \leq a/b \leq c_2 < \infty$ for some constants $c_1$ and $c_2$. If, for a fixed family of kernels $\mathcal{K} \in \mathcal{P}(\mathfrak{K})$, it is shown that $\mathbb{F}^* \subset \mathbb{F}(\mathcal{K})$ and

$$\sup_{F \in \mathbb{F}^*} \mathcal{U}_{\mathcal{K},F}(x) \asymp \inf_{\tilde{F}} \sup_{F \in \mathbb{F}^*} \mathcal{R}_r[\tilde{F}; F] \qquad \text{as } \varepsilon \to 0,$$

then the estimator $\hat{F}_{\mathcal{K}}$ is minimax on $\mathbb{F}^*$.

The minimax global results can be also derived from local inequalities of type (3). Indeed, suppose that we are interested in estimating $F$ with small $\mathbb{L}_r$-risk

$$\mathcal{R}_{\mathbb{L}_r}[\hat{F}; F] := \{\mathbb{E}_F \|\hat{F} - F\|_r^r\}^{1/r}, \qquad r > 0,$$

where $\|\cdot\|_r$ is the standard $\mathbb{L}_r$-norm on $\mathcal{D}_0$. Then, by the use of Fubini's theorem, we obtain from (3) that

$$\mathcal{R}_{\mathbb{L}_r}[\hat{F}_{\mathcal{K}}; F] \leq \|\mathcal{U}_{\mathcal{K},F}(\cdot)\|_r.$$

If, for a fixed family of kernels $\mathcal{K} \in \mathcal{P}(\mathfrak{K})$, one can prove that $\mathbb{F}^* \subset \mathbb{F}(\mathcal{K})$ and

$$\sup_{F \in \mathbb{F}^*} \|\mathcal{U}_{\mathcal{K},F}(\cdot)\|_r \asymp \inf_{\tilde{F}} \sup_{F \in \mathbb{F}^*} \mathcal{R}_{\mathbb{L}_r}[\tilde{F}; F] \qquad \text{as } \varepsilon \to 0,$$

then the corresponding estimator $\hat{F}_{\mathcal{K}}$ is minimax on $\mathbb{F}^*$ with respect to $\mathbb{L}_r$-risk. Local inequalities (3) are powerful tools for derivation of global minimax results in problems of estimating functions with inhomogeneous structure.

Local and global minimax adaptive results are obtained in a similar way. In the framework of the minimax adaptive approach, $F$ is assumed to belong to $\bigcup_{s \in S} \mathbb{F}_s^*$, where $\{\mathbb{F}_s^*, s \in S\}$ is a given collection of sets. The objective is to find an estimator $\hat{F}$ such that for every $s \in S$,

$$\sup_{F \in \mathbb{F}_s^*} \mathcal{R}_r[\hat{F}; F] \asymp \inf_{\tilde{F}} \sup_{F \in \mathbb{F}_s^*} \mathcal{R}_r[\tilde{F}; F] \qquad \text{as } \varepsilon \to 0. \tag{6}$$

If, for some $\mathcal{K} \in \mathcal{P}(\mathfrak{K})$, one can show that $\mathbb{F}_s^* \subset \mathbb{F}(\mathcal{K})$ for all $s \in S$ and

$$\sup_{F \in \mathbb{F}_s^*} \mathcal{U}_{\mathcal{K},F}(x) \asymp \inf_{\tilde{F}} \sup_{F \in \mathbb{F}_s^*} \mathcal{R}_r[\tilde{F}; F] \qquad \text{as } \varepsilon \to 0,$$

then the estimator $\hat{F}_{\mathcal{K}}$ is minimax adaptive for the collection $\{\mathbb{F}_s^*, s \in S\}$. Moreover, if $\mathbb{F}_s^* \subset \mathbb{F}(\mathcal{K}), \forall s \in S$ and

$$\sup_{F \in \mathbb{F}_s^*} \|\mathcal{U}_{\mathcal{K},F}(\cdot)\|_r \asymp \inf_{\tilde{F}} \sup_{F \in \mathbb{F}_s^*} \mathcal{R}_{\mathbb{L}_r}[\tilde{F}; F] \qquad \text{as } \varepsilon \to 0,$$

then $\hat{F}_{\mathcal{K}}$ is minimax adaptive for $\{\mathbb{F}_s^*, s \in S\}$ with respect to the $\mathbb{L}_r$-risk.

The rest of the paper is organized in the following way. In Section 2, we introduce notation and assumptions that are used throughout the paper and prove some preparatory results. In Section 3, we present our selection procedure, discuss its connections to other procedures and state the main result of this paper (Theorem 1). In Section 4, we apply the developed selection procedure to the aforementioned nonparametric estimation problems (i)–(iii). Section 5 contains the proof of Theorem 1. In Section 6, we prove all the results appearing in Section 4. Auxiliary results and proofs are collected in Appendix.

## 2. Preliminaries

We will use the following notation: $\|\cdot\|_p$ denotes the $\mathbb{L}_p(\mathcal{D}_0)$-norm, while $\|\cdot\|_{p,\infty}$ denotes the $\mathbb{L}_{p,\infty}(\mathbb{R}^d \times \mathcal{D}_0)$-norm,

$$\|G\|_{p,\infty} = \sup_{x \in \mathcal{D}_0} \left( \int_{\mathbb{R}^d} |G(t,x)|^p \, dt \right)^{1/p}, \qquad p \in [1,\infty].$$

We also write $|\cdot|_2$ for the Euclidean norm.

*Basic families of kernels.* Let $\Theta \subset \mathbb{R}^m$ be a compact set and consider a parameterized family of kernels $\mathcal{K}_\Theta = \{K_\mu, \mu \in \Theta\}$, where $K_\mu : \mathbb{R}^d \times \mathbb{R}^d \to \mathbb{R}$. Throughout the paper, we consider families of kernels $\mathcal{K}_\Theta$ satisfying the following conditions.

**(K0)** *Let $\mathcal{D}_1$ be an open interval in $\mathbb{R}^d$ such that $\mathcal{D}_0 \subset \mathcal{D}_1 \subset \mathcal{D}$. For all $\mu \in \Theta$, one has*

$$\begin{aligned} \operatorname{supp}(K_\mu(\cdot,y)) &\subseteq \mathcal{D}_1 \qquad \forall y \in \mathcal{D}_0, \\ \int_{\mathcal{D}} K_\mu(t,y) \, dt &= 1 \qquad \forall y \in \mathcal{D}_1. \end{aligned} \qquad (7)$$

*Moreover,*

$$\sigma(\mathcal{K}_\Theta) := \sup_{\mu \in \Theta} \|K_\mu\|_{2,\infty} < \infty, \qquad (8)$$

$$M(\mathcal{K}_\Theta) := \sup_{\mu \in \Theta} \|K_\mu\|_{1,\infty} < \infty. \qquad (9)$$

Note that (7) implies that $M(\mathcal{K}_\Theta) \geq 1$. Conditions (7)–(9) are standard in the context of kernel estimation.

In the construction of our selection rule, we use the auxiliary kernel collection $\mathcal{K}_{\Theta \times \Theta} = \{K_{\mu,\nu}, \mu, \nu \in \Theta\}$, $K_{\mu,\nu} : \mathbb{R}^d \times \mathbb{R}^d \to \mathbb{R}$, defined as

$$K_{\mu,\nu}(t,x) := \int_{\mathcal{D}_1} K_\mu(t,y) K_\nu(y,x) \, dy, \qquad t \in \mathcal{D}, x \in \mathcal{D}_0.$$

In what follows, we will assume that the following "commutativity property" is fulfilled for the kernels from $\mathcal{K}_\Theta$:

**(K1)**

$$K_{\mu,\nu} \equiv K_{\nu,\mu} \qquad \forall \mu, \nu \in \Theta. \tag{10}$$

*Remark 1.* Assumption (K1) is crucial for the construction of our selection procedure. Although this is a restriction on the family $\mathcal{K}_\Theta$, (10) is trivially fulfilled for kernels $K_\mu(t,x) = K_\mu(t-x)$ that correspond to standard kernel estimators.

The next statement establishes an important property of the kernel $K_{\mu,\nu}, \mu, \nu \in \Theta$. With any function $F$, we associate the quantities

$$B_{\mu,\nu}(x) = \int_\mathcal{D} K_{\mu,\nu}(t,x) F(t)\, \mathrm{d}t - F(x), \tag{11}$$

$$B_\nu(x) = \int_\mathcal{D} K_\nu(t,x) F(t)\, \mathrm{d}t - F(x), \qquad x \in \mathcal{D}_0. \tag{12}$$

**Proposition 1.** *Let (7) hold. Then for any $x \in \mathcal{D}_0$ and $F \in \mathbb{C}_b(\mathcal{D})$, one has*

$$B_{\mu,\nu}(x) - B_\nu(x) = \int_\mathcal{D} K_\nu(y,x) B_\mu(y)\, \mathrm{d}y. \tag{13}$$

The proof of the proposition is given in the Appendix.

*Remark 2.* Note that $K_{\nu,\mu}$ is a kernel for all $\mu, \nu \in \Theta$, that is, $\int_\mathcal{D} K_{\mu,\nu}(t,x)\, \mathrm{d}t = 1$, $\forall x \in \mathcal{D}_0$. This fact follows immediately from (13) if we put $F \equiv 1$.

*Auxiliary estimators and selection statistics.* With the collections $\mathcal{K}_\Theta$ and $\mathcal{K}_{\Theta \times \Theta}$, we associate the following families of linear estimators via (2):

$$\mathcal{F}(\mathcal{K}_\Theta) = \{\hat{F}_\mu = \hat{F}^{K_\mu}, \mu \in \Theta\};$$

$$\mathcal{F}(\mathcal{K}_{\Theta \times \Theta}) = \{\hat{F}_{\mu,\nu} = \hat{F}^{K_{\mu,\nu}}, \mu, \nu \in \Theta\}.$$

It is easily seen that

$$\hat{F}_\mu(x) - F(x) = B_\mu(x) + \varepsilon \xi_\mu(x),$$

$$\hat{F}_{\mu,\nu}(x) - F(x) = B_{\mu,\nu}(x) + \varepsilon \xi_{\mu,\nu}(x),$$

where

$$\xi_\mu(x) = \int K_\mu(t,x) W(\mathrm{d}t), \qquad \xi_{\mu,\nu}(x) = \int K_{\mu,\nu}(t,x) W(\mathrm{d}t).$$

Thus, the quantities $B_\mu(x)$ and $B_{\mu,\nu}(x)$ defined in (11)–(12) represent the bias of $\hat{F}_\mu(x)$ and $\hat{F}_{\mu,\nu}(x)$, respectively. In addition, we denote $\sigma_\mu^2(x) = \mathrm{var}\{\hat{F}_\mu(x)\} = \|K_\mu(\cdot, x)\|_2^2$.

Our selection procedure will be based on the statistics $\{\hat{F}_{\mu,\nu}(x) - \hat{F}_\nu(x), \mu, \nu \in \Theta\}$. It is clear that

$$\hat{F}_{\mu,\nu}(x) - \hat{F}_\nu(x) = B_{\mu,\nu}(x) - B_\nu(x) + \varepsilon[\xi_{\mu,\nu}(x) - \xi_\nu(x)], \tag{14}$$

where $\xi_{\mu,\nu}(x) - \xi_\nu(x)$ is a Gaussian zero-mean random variable with variance

$$\sigma^2_{\mu,\nu}(x) := \text{var}\{\hat{F}_{\mu,\nu}(x) - \hat{F}_\nu(x)\} = \|K_{\mu,\nu}(\cdot, x) - K_\nu(\cdot, x)\|_2^2.$$

Also, note that $\hat{F}_{\mu,\nu}(x) = \hat{F}_{\nu,\mu}(x)$, in view of (K1).

*Integrated bias and variance.* With any estimator $\hat{F}_\mu$, $\mu \in \Theta$, we associate the following two quantities:

$$\tilde{B}_\mu(x) := \sup_{\nu \in \Theta} |B_{\mu,\nu}(x) - B_\nu(x)| \vee |B_\mu(x)|; \tag{15}$$

$$\tilde{\sigma}_\mu(x) := \sup_{\nu \in \Theta} \int |K_\nu(y, x)| \sigma_\mu(y) \, dy \vee \sigma_\mu(x). \tag{16}$$

In words, $\tilde{B}_\mu$ is the maximum among the maximal integrated (with kernels $K_\nu$) bias and the bias of $\hat{F}_\mu$, while $\tilde{\sigma}_\mu$ is the maximum among the maximal integrated standard deviation of $\hat{F}_\mu$ and standard deviation of $\hat{F}_\mu$. In what follows, with slight abuse of terminology, we will refer to $\tilde{B}_\mu(x)$ and $\tilde{\sigma}_\mu^2(x)$ as the *integrated bias* of $\hat{F}_\mu$ and the *integrated variance* of $\hat{F}_\mu$, respectively.

It follows from (13) and (9) that

$$\tilde{B}_\mu(x) \leq M(\mathcal{K}_\Theta) \sup_y |B_\mu(y)|, \qquad \tilde{\sigma}_\mu(x) \leq M(\mathcal{K}_\Theta) \sup_y \sigma_\mu(y). \tag{17}$$

We also have the following upper bound on $\sigma_{\mu,\nu}(x)$ in terms of $\tilde{\sigma}_\mu(x)$ and $\tilde{\sigma}_\nu(x)$: for all $\mu, \nu \in \Theta$,

$$\begin{aligned}\sigma_{\mu,\nu}(x) &\leq \|K_{\mu,\nu}(\cdot, x)\|_2 + \|K_\nu(\cdot, x)\|_2 \\ &\leq \tilde{\sigma}_\mu(x) + \sigma_\nu(x) \leq \tilde{\sigma}_\mu(x) + \tilde{\sigma}_\nu(x).\end{aligned} \tag{18}$$

Here, we have used the triangle inequality and the Minkowski inequality for integrals.

In what follows, point $x$ is fixed. So, in our notation, we will not indicate dependence on $x$ when this does not lead to confusion.

## 3. Selection procedure and main result

In this section, we introduce our selection rule.

### 3.1. Majorant

We begin with the definition of the *majorant*, the main ingredient of our construction.

Let $\Sigma_\Theta := \{\tilde{\sigma}_\mu : \mu \in \Theta\} \subset \mathbb{R}_+$ and define $\sigma_{\min} := \inf \Sigma_\Theta$, $\sigma_{\max} := \sup \Sigma_\Theta$. Thus, $\Sigma_\Theta$ is the image of $\Theta$ under the mapping $\mu \mapsto \tilde{\sigma}_\mu$, where $\tilde{\sigma}_\mu$ is defined in (16). Let

$$e_{\mathcal{K}_\Theta}(\sigma) := \sup_{\mu \in \Theta} \mathbb{E} \sup_{\nu : \tilde{\sigma}_\nu \leq \sigma} |\xi_{\mu,\nu} - \xi_\nu|, \qquad \sigma \in \Sigma_\Theta. \tag{19}$$

*Remark 3.* By definition, the function $e_{\mathcal{K}_\Theta}(\cdot)$ is non-decreasing on $\Sigma_\Theta$. For any given $\sigma \in \Sigma_\Theta$, $e_{\mathcal{K}_\Theta}(\sigma)$ is the maximal (over $\mu \in \Theta$) expectation of supremum of the Gaussian zero-mean random process $\{\xi_{\mu,\nu} - \xi_\nu\}$ with the index set $\{\nu : \tilde{\sigma}_\nu \leq \sigma\} \subseteq \Theta$. The covariance structure of this process is completely determined by the family of kernels $\mathcal{K}_\Theta$. Thus, the function $e_{\mathcal{K}_\Theta}(\cdot)$ can be computed, for example, using Monte Carlo simulations. Alternatively, useful analytical bounds on $e_{\mathcal{K}_\Theta}(\cdot)$ can be derived from the theory of Gaussian processes.

**(E)** *Let $e(\sigma)$ be a continuous non-decreasing function on $\Sigma_\Theta$ such that*

  (i) $e(\sigma) \geq e_{\mathcal{K}_\Theta}(\sigma)$, $\forall \sigma \in \Sigma_\Theta$,
  (ii) *there exist absolute constants $1 < c_e \leq C_e$ such that*

$$c_e \leq \frac{e(2\sigma)}{e(\sigma)} \leq C_e \qquad \forall \sigma \in \Sigma_\Theta. \tag{20}$$

*Remark 4.* The function $e(\cdot)$ is an upper bound on $e_{\mathcal{K}_\Theta}(\cdot)$. Such a bound can be derived from general inequalities on suprema of Gaussian processes. Condition (20) holds, for example, if $e(\sigma) = c\sigma L(\sigma)$, where $c$ is a constant and $L(\sigma)$ is a slowly varying function. In fact, for our purposes, it is sufficient to require that inequalities in (20) hold for the ratio $e(a\sigma)/e(\sigma)$ for some $a > 1$.

We are now in a position to define the *majorant*:

$$Q(\sigma) := \varkappa_0 e(\sigma) + \sigma \sqrt{1 + \varkappa_1 \ln \frac{\sigma}{\sigma_{\min}}}, \qquad \sigma \in \Sigma_\Theta, \tag{21}$$

where $\varkappa_0 = 2C_e$ and $\varkappa_1 = 128r(1 \vee \ln C_e / \ln 2)$.

*Remark 5.* Loosely speaking, the majorant uniformly bounds from above the random process $\xi_{\mu,\nu} - \xi_\nu$, $\mu, \nu \in \Theta$, with prescribed probability. The function $Q$ consists of two terms. The first term bounds the expectation of the supremum of a zero-mean Gaussian random process, while the second term controls the deviation of this supremum from its expectation. In fact, the first term characterizes "massiveness" of the subset of estimators from $\mathcal{F}(\mathcal{K}_\Theta)$ with variance less than a prescribed level. The second term involves a logarithm of the ratio of estimator variances in the family. It can be regarded as a price to be paid for considering families of estimators with different variances.

### 3.2. Selection rule

We are now in a position to define our selection rule.

For any $\mu \in \Theta$, let

$$\hat{R}_\mu := \sup_{\nu: \tilde{\sigma}_\nu \geq \tilde{\sigma}_\mu} \{|\hat{F}_{\mu,\nu} - \hat{F}_\nu| - \tfrac{1}{2}\varepsilon Q(\tilde{\sigma}_\nu)\}. \qquad (22)$$

Let $\delta = \tfrac{1}{4}\varepsilon Q(\sigma_{\min})$ and let $\hat{\mu} \in \Theta$ be such that

$$\hat{R}_{\hat{\mu}} + \varepsilon Q(\tilde{\sigma}_{\hat{\mu}}) \leq \inf_{\mu \in \Theta}\{\hat{R}_\mu + \varepsilon Q(\tilde{\sigma}_\mu)\} + \delta. \qquad (23)$$

We then define

$$\hat{F} = \hat{F}_{\hat{\mu}}. \qquad (24)$$

Several remarks on the above definition are in order. First, observe that $\hat{R}_\mu$ may be negative; however, by definition,

$$\hat{R}_\mu \geq -\tfrac{1}{2}\varepsilon Q(\tilde{\sigma}_\mu) \qquad \forall \mu \in \Theta, \qquad (25)$$

so that $\hat{R}_\mu + \varepsilon Q(\tilde{\sigma}_\mu)$ is always positive. Second, in order to ensure that there exists a measurable choice of $\hat{\mu}$ satisfying (23), one needs to impose additional conditions on the family of kernels $\mathcal{K}_\Theta$. The next assumption provides such conditions.

**(K2)** *There exist positive constants $\bar{L}$ and $\gamma \in (0,1]$ such that*

$$\sup_{\mu,\mu' \in \Theta} \frac{\|\widetilde{K}_\mu - \widetilde{K}_{\mu'}\|_{2,\infty}}{|\mu - \mu'|_2^\gamma} \leq \bar{L}, \qquad (26)$$

$$\sup_{\mu,\mu' \in \Theta} \frac{\sup_x |1 - \|K_\mu(\cdot,x)\|_2/\|K_{\mu'}(\cdot,x)\|_2|}{|\mu - \mu'|_2^\gamma} \leq \bar{L}, \qquad (27)$$

*where $\widetilde{K}_\mu(\cdot,x) = K_\mu(\cdot,x)/\|K_\mu(\cdot,x)\|_2$, $\forall \mu \in \Theta$.*

In the proof of Theorem 1, we show that (K0)–(K2), and boundedness and continuity of $F$ imply that there exists a measurable choice of $\hat{\mu} \in \Theta$ such that (23) holds. Thus, our selection rule is well defined.

### 3.3. Discussion

In this section, we explain the main idea underlying the construction of our selection scheme and discuss connections to other procedures in the literature.

The pointwise selection procedures were developed by Lepski (1990, 1991), Lepski *et al.* (1997), Lepski and Spokoiny (1997) and Kerkyacharian *et al.* (2001). In those papers, the procedures are two-staged: first, a collection of *admissible* estimators is constructed using a "bias–variance" comparison scheme; second, among admissible estimators, an estimator with minimal variance is selected. The procedure in Lepski (1990, 1991) (and its refinements in Lepski *et al.* (1997), Lepski and Spokoiny (1997)) selects from the collection $\mathcal{F}(\mathcal{K}_1)$ of one-dimensional kernel estimators (see Example 1 in Section 1) discretized in an appropriate way. In our notation, it reads as follows:

select the estimator with maximal bandwidth $\mu \in [h_{\min}, h_{\max}]^{\#}$ such that

$$|\hat{F}_\mu - \hat{F}_\nu| \leq T(\mu, \nu) \qquad \forall \nu \in [h_{\min}, h_{\max}]^{\#} : \sigma_\nu \geq \sigma_\mu, \tag{28}$$

where $A^{\#}$ stands for a discretization of a set $A$ and $T(\mu, \nu)$ is a certain threshold.

Here, the set of admissible estimators contains all estimators $\hat{F}_\mu$, $\mu \in [h_{\max}, h_{\min}]^{\#}$ satisfying (28) and at the selection stage, the estimator with minimal variance (maximal bandwidth) is chosen. This scheme exploits monotonicity properties of the bias and variance with respect to the bandwidth which, in general, do not hold in the multidimensional case.

A generalization of (28) to the multidimensional case was developed in Kerkyacharian *et al.* (2001). Their procedure is designed for selection from the properly discretized collection $\mathcal{F}(\mathcal{K}_\mathcal{H})$ (see Example 2, Section 1) and can be represented as follows:

For $\mu = (\mu_1, \ldots, \mu_d) \in \mathcal{H}^{\#}$ and $\nu = (\nu_1, \ldots, \nu_d) \in \mathcal{H}^{\#}$, define $\mu \vee \nu = (\mu_1 \vee \nu_1, \ldots, \mu_d \vee \nu_d)$ and consider the auxiliary estimator $\hat{F}_{\mu,\nu} \equiv \hat{F}_{\mu \vee \nu}$. The estimator $\hat{F}_\mu$, $\mu \in \mathcal{H}^{\#}$, is called *admissible* if

$$|\hat{F}_{\mu,\nu} - \hat{F}_\nu| \leq T(\nu) \qquad \forall \nu \in \mathcal{H}^{\#} : \sigma_\nu \geq \sigma_\mu, \tag{29}$$

where $T(\nu)$ is an appropriate threshold. Note that (29) can be rewritten as

$$\sup_{\nu \in \mathcal{H}^{\#} : \sigma_\nu \geq \sigma_\mu} [|\hat{F}_{\mu,\nu} - \hat{F}_\nu| - T(\nu)] \leq 0. \tag{30}$$

At the selection stage, we choose the admissible estimator with minimal variance.

Note that the scheme (29) involves an auxiliary estimator $\hat{F}_{\mu,\nu}$ and its construction can only be used for selection from the collection $\mathcal{F}(\mathcal{K}_\mathcal{H})$. Specifically, the procedure (29) cannot be applied for selection from the collection of kernel estimators $\mathcal{F}(\mathcal{K}_{SI})$ (see Example 3, Section 1) corresponding to the single index model.

Our selection procedure (22)–(24) also uses an auxiliary estimator $\hat{F}_{\mu,\nu}$, but, in contrast to (29), the construction of $\hat{F}_{\mu,\nu}$ is universal and fits a wide variety of kernel collections. In addition, instead of pairwise comparisons with a threshold (as in (28) and (29)), we define the majorant function and use direct minimization. Our rule (23) is very much in the spirit of (30). Indeed, the procedure of Kerkyacharian *et al.* (2001) minimizes $\sigma_\mu$ subject to constraint (30), while in (23), we minimize, with respect to $\mu$, the expression $\sup_{\nu : \sigma_\nu \geq \sigma_\mu} [|\hat{F}_{\mu,\nu} - \hat{F}_\mu| - \frac{1}{2}T(\nu)] + T(\mu)$ and $T(\mu)$ is "roughly" proportional to $\sigma_\mu$.

Summing up, the proposed selection method differs from other pointwise selection procedures in: (a) construction of the auxiliary estimators $\hat{F}_{\mu,\nu}$; (b) selection by direct minimization. These features enable a wide variety of kernel collections to be treated in a unified way and the discretization of the parameter space $\Theta$ to be avoided.

### 3.4. Upper bound

In order to present an upper bound on the risk of the proposed estimator, we need the following definition.

For any function $F \in \mathbb{C}_b(\mathcal{D})$ and given collection $\mathcal{K}_\Theta$, define

$$\Theta_F(\mathcal{K}_\Theta) := \{\mu \in \Theta : \forall \sigma \geq \tilde{\sigma}_\mu, \sigma \in \Sigma_\Theta \; \exists \theta \in \Theta \text{ such that } \tilde{\sigma}_\theta = \sigma \text{ and } \tilde{B}_\theta \leq \tfrac{1}{4}\varepsilon Q(\tilde{\sigma}_\theta)\}.$$

In what follows, we will consider functions $F$ for which $\Theta_F(\mathcal{K}_\Theta)$ is non-empty. This condition is closely related to the existence of estimators in $\mathcal{F}(\mathcal{K}_\Theta)$ realizing the bias-variance trade-off.

*Remark 6.* Clearly, $\Theta_F(\mathcal{K}_\Theta)$ is non-empty for any constant function $F$ since, by (13) and (15), $\tilde{B}_\theta \equiv 0$ for all $\theta \in \Theta$. For the same reason, if $K_\theta$ is orthogonal to all polynomials of degree $\leq l$, then $\Theta_F(\mathcal{K}_\Theta)$ is non-empty for any $F$ which is a polynomial of degree $\leq l$. In general, the size of the set of functions $F$ for which $\Theta_F(\mathcal{K}_\Theta)$ is non-empty is completely determined by the family $\mathcal{K}_\Theta$. For example, if $\mathcal{F}(\mathcal{K}_\Theta)$ is the family of standard kernel estimators with a bounded kernel $K_\theta$ and bandwidth $\theta = (h_1, \ldots, h_d) \in [\varepsilon^2, 1/2]^d$, then $\Theta_F(\mathcal{K}_\Theta)$ is non-empty for any $F \in \mathbb{C}_b(\mathcal{D})$.

Finally, we put

$$\mu^* = \arg\inf_{\mu \in \Theta_F(\mathcal{K}_\Theta)} \tilde{\sigma}_\mu. \tag{31}$$

**Theorem 1.** *Suppose that assumptions* (K0)–(K2) *and* (E) *hold. Then, for any* $F \in \mathbb{C}_b(\mathcal{D})$ *such that* $\Theta_F(\mathcal{K}_\Theta) \neq \varnothing$, *and for all* $\varepsilon$ *small enough, one has*

$$\mathcal{R}_r[\hat{F}; F] \leq C\varepsilon Q(\tilde{\sigma}_{\mu^*}),$$

*where* $C$ *is a numerical constant depending only on* $r$, $c_e$ *and* $C_e$.

## 4. Applications

In this section, we show how the upper bound of Theorem 1 can be used for the derivation of minimax and adaptive minimax results. In particular, in Sections 4.2–4.4, we consider three particular problems:

- pointwise adaptive estimation in the single index model;
- pointwise minimax estimation over a union of anisotropic Hölder classes;
- global minimax estimation over isotropic Besov classes.

Our goal here is to show how a careful choice of the family of kernels leads to estimators with optimal statistical properties. Note that in each particular case, the estimators are different, although all of them are obtained by the same computational routine presented in Section 3.

In Section 4.5, we demonstrate that the choice of a rather huge kernel collection allows a *single* estimator to be constructed which is simultaneously optimal (up to a log-factor) for these three entirely different problems.

## 4.1. General kernel collection

Let $G:\mathbb{R}^d \to \mathbb{R}$ be a function supported on $[-1/2, 1/2]^d$ and satisfying the conditions

$$\int G(t)\,\mathrm{d}t = 1, \qquad \int t^{\boldsymbol{r}} G(t)\,\mathrm{d}t = 0 \quad \forall |\boldsymbol{r}| = 1,\ldots,l, \qquad \sup_t |\nabla G(t)|_2 \leq M, \quad (32)$$

where $\boldsymbol{r} = (r_1,\ldots,r_d)$, $r_i \geq 0$, $|\boldsymbol{r}| = r_1 + \cdots + r_d$ and $t^{\boldsymbol{r}} = t_1^{r_1}\cdots t_d^{r_d}$.

Let $\mathcal{E}$ denote the set of $d\times d$ orthogonal matrices and let $\mathcal{H} = [h_{\min}, h_{\max}]^d$, where $0 < h_{\min} \leq h_{\max} \leq 1/2$ are given real numbers.

Define, for all $h \in \mathcal{H}$ and all $E \in \mathcal{E}$,

$$G_h(t) = \left[\prod_{i=1}^d h_i^{-1}\right] G\left(\frac{t_1}{h_1},\ldots,\frac{t_d}{h_d}\right), \qquad G_{h,E}(t,x) = G_h(E^T[t-x]) \quad (33)$$

and consider the following collection of kernels:

$$\mathcal{K}_{\mathcal{H},\mathcal{E}} = \{G_{h,E}, (h,E) \in \mathcal{H}\times\mathcal{E}\}. \quad (34)$$

***Remark 7.***

1. For the family $\mathcal{K}_{\mathcal{H},\mathcal{E}}$, we have

$$\sigma_{h,E}(x) = \sigma_{h,E} = \|G\|_2 \prod_{i=1}^d h_i^{-1/2} \qquad \forall x \in \mathcal{D}_0,\ \forall E \in \mathcal{E}$$

and, therefore, $\tilde{\sigma}_{h,E}(x) = \|G\|_1 \sigma_{h,E}$, $\forall x \in \mathcal{D}_0$,

$$\sigma_{\min} = \|G\|_1 \|G\|_2 \, h_{\max}^{-d/2}, \qquad \sigma_{\max} = \|G\|_1 \|G\|_2 \, h_{\min}^{-d/2}.$$

2. Assumptions (K0)–(K2) are fulfilled for the family $\mathcal{K}_{\mathcal{H},\mathcal{E}}$. Indeed, (K0) holds trivially; here, $M(\mathcal{K}_{\mathcal{H},\mathcal{E}}) = \|G\|_1$. Assumption (K2) is fulfilled because $\mathcal{K}_{\mathcal{H},\mathcal{E}}$ consists of convolution kernels. Boundedness of the gradient of $G$ in (32), along with (K0) and (K1), implies (K2) (see Lemmas 1 and 2 in Section 6.1).

In order to construct estimators in the aforementioned problems, we will consider families corresponding to different subsets of $\mathcal{K}_{\mathcal{H},\mathcal{E}}$. The family of estimators $\mathcal{F}(\mathcal{K}_{\mathcal{H},\mathcal{E}})$ will be considered in Section 4.5.

### 4.2. Pointwise adaptive estimation in the single index model

Consider the model (1) with $F(t) = f(\omega^T t)$, where $f : \mathbb{R} \to \mathbb{R}$ is an unknown function from the Hölder ball $\mathbb{H}_1(\alpha, L)$ with unknown parameters $\alpha > 0$ and $L > 0$, and $\omega \in \mathbb{S}^{d-1}$ is an unknown direction vector. We refer to this model as the *single index model*.

**Definition 1.** *We say that function $F$ belongs to the functional class $\mathbb{F}_{SI}(\alpha, L)$ if there exists a direction vector $\omega \in \mathbb{S}^{d-1}$, parameters $\alpha > 0$, $L > 0$ and a univariate function $f \in \mathbb{H}_1(\alpha, L)$ such that $F(t) = f(\omega^T t)$.*

Define $\mathcal{H}_1 = \{h \in \mathcal{H} : h = (h_1, h_{\max}, \ldots, h_{\max})\}$, $\Theta_{SI} = \mathcal{H}_1 \times \mathcal{E}$ and consider the following subset of $\mathcal{K}_{\mathcal{H},\mathcal{E}}$:

$$\mathcal{K}_{SI} = \{G_{h,E} : (h, E) \in \Theta_{SI}\}. \tag{35}$$

The corresponding family of estimators is given by

$$\mathcal{F}(\mathcal{K}_{SI}) = \left\{ \hat{F}_{h,E}(x) = \int G_{h,E}(t, x) Y(dt), (h, E) \in \Theta_{SI} \right\}.$$

**Remark 8.** In view of Remark 7, we have

$$\tilde{\sigma}_{h,E} = \|G\|_1 \sigma_{h,E} = \|G\|_1 \|G\|_2 \, h_{\max}^{-(d-1)/2} [1/\sqrt{h_1}]; \tag{36}$$

$$\sigma_{\min} = h_{\max}^{-d/2} \|G\|_1 \|G\|_2, \qquad \sigma_{\max} = h_{\max}^{-(d-1)/2} \|G\|_1 \|G\|_2 [1/\sqrt{h_{\min}}]. \tag{37}$$

Note that $\sigma_{h,E}$ does not depend on $E$.

Let $e(\sigma) = C_0 \sigma \sqrt{\ln \sigma}$, where $C_0$ is a numerical constant depending only on $d$ and $G$. It is shown in Lemma 3 of Section 6 that $e(\sigma) \geq e_{\mathcal{K}_{SI}}(\sigma)$ for all $\sigma \in \Sigma_{\Theta_{SI}}$. The majorant $Q$ is given by

$$Q(\sigma) = \sigma[\varkappa_0 C_0 \sqrt{\ln \sigma} + \sqrt{1 + \varkappa_1 \ln(\sigma/\sigma_{\min})}]. \tag{38}$$

Note that assumption (E) is trivially fulfilled with $c_e = 2$ and

$$C_e = 2(1 + \sqrt{\ln 2 / \ln \sigma_{\min}}) \leq 2(1 + \sqrt{2/d}).$$

Let $\hat{F}_{SI}$ be the estimator derived from the collection $\mathcal{F}(\mathcal{K}_{SI})$, in accordance with our general selection rule, with the majorant (38), where $\varkappa_0 = 4(1 + \sqrt{\ln 2 / \ln \sigma_{\min}})$ and $\varkappa_1 = 320r$.

**Theorem 2.** *Fix some $0 < \alpha_{\max} < \infty$, let $h_{\min} = \varepsilon^2$, $h_{\max} = \varepsilon^{2/(2\alpha_{\max}+1)}$ and assume that (32) holds with $l \geq \lfloor \alpha_{\max} \rfloor$.*

*Then, for any $0 < \alpha \leq \alpha_{\max} \leq l$, $L > 0$ and $\varepsilon$ small enough, one has*

$$\sup_{F \in \mathbb{F}_{SI}(\alpha, L)} \mathcal{R}_r[\hat{F}_{SI}; F] \leq C L^{1/(2\alpha+1)} \left(\varepsilon \sqrt{\ln \frac{1}{\varepsilon}}\right)^{2\alpha/(2\alpha+1)}, \tag{39}$$

*where $C$ depends only on $\alpha$, $G$, $d$ and $r$.*

**Remark 9.** If the parameters $\alpha$ and $L$ of the class $\mathbb{H}_1(\alpha, L)$ are known and the direction vector $\omega$ is unknown, then we consider the following subset of $\mathcal{K}_{SI}$:

$$\mathcal{K}'_{SI} = \{G_{h,E}, h = h^*, E \in \mathcal{E}\},$$

where $h^* = (h_1^*, h_{\max}, \ldots, h_{\max})$, $h_1^* = L^{-2/(2\alpha+1)}[\varepsilon\sqrt{\ln 1/\varepsilon}]^{2/(2\alpha+1)}$. Under these circumstances,

$$\sigma_{h^*, E} = h_{\max}^{-(d-1)/2}(h_1^*)^{-1/2}\|G\|_2$$

does not depend on $E$ (see also Remark 8) and therefore

$$\sigma_{\min} = \sigma_{\max} = \sigma^* := \|G\|_1 \sigma_{h^*, E} = h_{\max}^{-(d-1)/2}(h_1^*)^{-1/2}\|G\|_2 \|G\|_1.$$

The corresponding majorant is given by $Q(\sigma^*) = \varkappa_0 C_0 \sigma^* \sqrt{\ln \sigma^*} + \sigma^*$ so that the first term is dominating (all estimators in $\mathcal{F}(\mathcal{K}'_{SI})$ have the same variance). The resulting selected estimator for this family will then satisfy the same upper bound of Theorem 2. One can prove a lower bound that shows that even if $\alpha$ and $L$ are known, the rate of convergence on the right-hand side of (39) cannot be improved.

### 4.3. Pointwise minimax estimation over a union of anisotropic Hölder classes

We start with the definition of the anisotropic Hölder class of functions.

**Definition 2.** *Let $\boldsymbol{\alpha} = (\alpha_1, \ldots, \alpha_d)$, $\alpha_i > 0$ and $L > 0$. We say that $f:[-1/2, 1/2]^d \to \mathbb{R}$ belongs to the anisotropic Hölder class $\mathbb{H}_d(\boldsymbol{\alpha}, L)$ if for all $i = 1, \ldots, d$ and all $t \in [-1/2, 1/2]^d$,*

$$|D_i^m f(t)| \leq L \qquad \forall m = 1, \ldots, \lfloor \alpha_i \rfloor$$

*and*

$$|D_i^{\lfloor \alpha_i \rfloor} f(t_1, \ldots, t_i + z, \ldots, t_d) - D_i^{\lfloor \alpha_i \rfloor} f(t_1, \ldots, t_i, \ldots, t_d)| \leq L|z|^{\alpha_i - \lfloor \alpha_i \rfloor} \qquad \forall z \in \mathbb{R},$$

*where $D_i^m f$ denotes the mth order partial derivative of $f$ with respect to the variable $t_i$ and $\lfloor \alpha_i \rfloor$ is the largest integer strictly less than $\alpha_i$.*

Fix $\gamma > 0$ and introduce the functional class

$$\mathbb{F}_{AH}(\gamma, L) = \bigcup_{\boldsymbol{\alpha} \in \mathcal{A}_\gamma} \mathbb{H}_d(\boldsymbol{\alpha}, L), \qquad \text{where } \mathcal{A}_\gamma = \left\{ \boldsymbol{\alpha} : \sum_{i=1}^{d} 1/\alpha_i = 1/\gamma, \alpha_i > 0, i = \overline{1, d} \right\}.$$

***Remark 10.*** It is well known (see, e.g., Kerkyacharian *et al.* (2001) and Bertin (2004)) that for any $\boldsymbol{\alpha} \in \mathcal{A}_\gamma$, the minimax rate of convergence on $\mathbb{H}_d(\boldsymbol{\alpha}, L)$ is given by $\varepsilon^{2\gamma/(2\gamma+1)}$. Thus, $\mathbb{F}_{AH}(\gamma, L)$ is the union of functional classes with prescribed accuracy of estimation. Klutchnikoff (2005) showed that the rate $\varepsilon^{2\gamma/(2\gamma+1)}$ is not achievable on $\mathbb{F}_{AH}(\gamma, L)$ and proved that the minimax rate of convergence on $\mathbb{F}_{AH}(\gamma, L)$ is given by

$$\varphi_\varepsilon := [\varepsilon \sqrt{\ln \ln(1/\varepsilon)}]^{2\gamma/(2\gamma+1)}.$$

In this section, we show that the application of our general selection rule with a specific choice of the kernel collection $\mathcal{K}_{AH} \subset \mathcal{K}_{\mathcal{H},\mathcal{E}}$ leads to the minimax estimator on $\mathbb{F}_{AH}(\gamma, L)$.

Define the set of bandwidths $\mathcal{H}_\gamma \subset \mathcal{H}$

$$\mathcal{H}_\gamma := \left\{ h \in [h_{\min}, h_{\max}]^d : \prod_{i=1}^{d} h_i^\gamma = \varphi_\varepsilon \right\} \tag{40}$$

and consider the following subset of the family of kernels $\mathcal{K}_{\mathcal{H},\mathcal{E}}$:

$$\mathcal{K}_{AH} = \{G_{h,E} : (h, E) \in \Theta_{AH} := \mathcal{H}_\gamma \times \{I_d\}\}, \tag{41}$$

where $I_d$ is the $d \times d$ identity matrix.

The corresponding family of estimators is given by

$$\mathcal{F}(\mathcal{K}_{AH}) = \left\{ \hat{F}_h(x) = \int G_h(t - x) Y(\mathrm{d}t), h \in \mathcal{H}_\gamma \right\}.$$

For all $h \in \mathcal{H}_\gamma$, we have

$$\tilde{\sigma}_h = \|G\|_1 \sigma_h = \|G\|_1 \|G\|_2 \prod_{i=1}^{d} h_i^{-1/2} = \|G\|_1 \|G\|_2 [\varepsilon \sqrt{\ln \ln(1/\varepsilon)}]^{-1/(2\gamma+1)}.$$

Thus, the set $\Sigma_{\Theta_{AH}}$ consists of the single point $\|G\|_1 \|G\|_2 [\varepsilon \sqrt{\ln \ln(1/\varepsilon)}]^{-1/(2\gamma+1)}$.

Let $e(\sigma) = C_1 \sigma \sqrt{\ln \ln(h_{\max}/h_{\min})}$, where $C_1$ is a numerical constant depending only on $d$ and $G$. Lemma 3 of Section 6 shows that $e(\sigma)$ is an upper bound on $e_{\mathcal{K}_{AH}}(\sigma)$. Note that assumption (E) is trivially fulfilled with $c_e = C_e = 2$ and the majorant in our procedure can be taken as follows:

$$Q(\sigma) = \sigma[1 + 4C_1 \sqrt{\ln \ln(h_{\max}/h_{\min})}]. \tag{42}$$

Let $\hat{F}_{AH}$ be the estimator derived from the collection $\mathcal{F}(\mathcal{K}_{AH})$, in accordance with our general selection rule, with the majorant (42).

**Theorem 3.** *Fix $0 < \alpha_{\max} < \infty$. Let $h_{\min} = \varepsilon^2$, $h_{\max} = 1/2$ and assume that (32) holds with $l \geq \lfloor \alpha_{\max} \rfloor$. Then, for any $\boldsymbol{\alpha} \in \mathcal{A}_\gamma \cap (0, \alpha_{\max}]^d$,*

$$\sup_{F \in \mathbb{H}_d(\boldsymbol{\alpha}, L)} \mathcal{R}_r[\hat{F}_{AH}; F] \leq CL^{1/(2\gamma+1)} \left( \varepsilon \sqrt{\ln \ln \frac{1}{\varepsilon}} \right)^{2\gamma/(2\gamma+1)},$$

*where $C$ depends only on $G$, $d$, $r$ and $\gamma$.*

### 4.4. Global minimax estimation over isotropic Besov classes

We begin with the definition of the isotropic Besov class of functions on $\mathcal{D}_0$.

For all $x \in \mathcal{D}_0$ and $a \in \mathbb{R}^d$ such $x + a \in \mathcal{D}$, define

$$\Delta_a^1 F(x) = F(x+a) - F(x).$$

For any integer $l \geq 2$, let $\Delta_a^l F(x)$ denote the $(l-1)$-fold iteration of the operator $\Delta_a^1 F(x)$.

**Definition 3.** *Let $s > 0, p \in [1, \infty)$ and $L > 0$ be given constants. Let $\mathbb{B}_{p,\infty}^s(d, L)$ denote the set of all functions satisfying*

$$\sup_a |a|_2^{-s} \|\Delta_a^{\lfloor s \rfloor + 2} F(\cdot)\|_p \leq L,$$

*where $\lfloor s \rfloor$ is the largest integer strictly less than $s$. We call $\mathbb{B}_{p,\infty}^s(d, L)$ the* isotropic Besov class *of functions.*

The considered classes were first introduced in approximation theory by Nikolskii (1975). They represent a particular case of the *Besov classes* $\mathbb{B}_{p,q}^s(d, L)$ with $q = \infty$ which appear more often in the statistical literature. More general *anisotropic* Besov functional classes were considered in Kerkyacharian *et al.* (2001).

On the class $\mathbb{B}_{p,\infty}^s(d, L)$, we introduce the *maximal risk*

$$\mathcal{R}_{\mathbb{L}_r}(\tilde{F}) = \sup_{F \in \mathbb{B}_{p,\infty}^s(d,L)} \{\mathbb{E}_F \|\tilde{F} - F\|_r^r\}^{1/r}, \qquad r \in [1, \infty),$$

where $\tilde{F}$ is an estimator of $F$. It is well known (Delyon and Juditsky (1996)) that

$$\varphi_\varepsilon = \begin{cases} \varepsilon^{s/(s+d/2)}, & \text{if } sp > \dfrac{d(r-p)}{2}, \\ [\varepsilon \sqrt{\ln 1/\varepsilon}]^{s/(s+d/2)} [\ln 1/\varepsilon]^{1/r}, & \text{if } sp = \dfrac{d(r-p)}{2}, \\ [\varepsilon \sqrt{\ln 1/\varepsilon}]^{s-d(1/p-1/r)/(s-d(1/p-1/2))}, & \text{if } sp < \dfrac{d(r-p)}{2} \end{cases}$$

is the minimax rate of convergence on $\mathbb{B}_{p,\infty}^s(d,L)$ if $sp \neq \frac{d(r-p)}{2}$ and differs from minimax rate of convergence by $\ln(1/\varepsilon)$-factor if $sp = \frac{d(r-p)}{2}$.

In this section, we present the estimator which attains the rate $\varphi_\varepsilon$ on $\mathbb{B}_{p,\infty}^s(d,L)$. As before, this estimator is the output of our general selection procedure.

Let

$$G(t) = G^*(t) := \sum_{j=1}^{\lfloor s \rfloor + 2} (-1)^{j+1} \binom{\lfloor s \rfloor + 2}{j} \frac{1}{j^d} g\left(\frac{t}{j}\right), \qquad t \in \mathbb{R}^d,$$

where $g : \mathbb{R}^d \to \mathbb{R}$ is a bounded, compactly supported function with $\int g = 1$. It is easily seen that the function $G^*$ satisfies assumption (32).

Consider the following subset of $\mathcal{K}_{\mathcal{H},\mathcal{E}}$:

$$\mathcal{K}_B = \{G_{h,E}^* : (h,E) \in \Theta_B := \mathcal{H}_B \times \{I_d\}\},$$

where $\mathcal{H} \supset \mathcal{H}_B := \{h = (h_1, \ldots, h_d) \in \mathcal{H} : h_i = h_j, i, j = \overline{1,d}\}$. Note that the family $\mathcal{K}_B$ consists of *isotropic kernels* having the same bandwidth in each direction. The corresponding family of estimators is given by

$$\mathcal{F}(\mathcal{K}_B) = \left\{\hat{F}_h : \hat{F}_h(x) = \int G_h^*(t-x) Y(\mathrm{d}t), h \in \mathcal{H}_B\right\}.$$

Let $\hat{F}_B$ be the estimator derived from the collection $\mathcal{F}(\mathcal{K}_B)$ in accordance with our general selection rule, where the majorant $Q$ is given by

$$Q(\sigma) = C(s,g)\sigma\sqrt{1 + \varkappa_1 \ln(\sigma/\sigma_{\min})} =: \mathcal{C}_1 \sigma Q^*(\sigma/\sigma_{\min}).$$

Here, $Q^*(z) = z\sqrt{1 + \ln z}, z \geq 1$ and $\mathcal{C}_1 = C(s,g,d,r)$ is the numerical constant.

**Theorem 4.** *Suppose that $s > d/p$ and choose $h_{\min} = \varepsilon^2$ and*

$$h_{\max} = \begin{cases} \varepsilon^{2/(2s+d)}, & \text{if } sp > \dfrac{d(r-p)}{2}, \\ 1/2, & \text{if } sp \leq \dfrac{d(r-p)}{2}. \end{cases}$$

*Then, for all $\varepsilon > 0$ small enough,*

$$\mathcal{R}_{\mathbb{L}_r}(\hat{F}_B) \leq C(d,s,p,r,g)\varphi_\varepsilon^r,$$

*where $C(d,s,p,r,g) > 0$ is a numerical constant.*

*Remark 11.* The result described in Theorem 4 was first obtained by Delyon and Juditsky (1996) using wavelet techniques. Lepski *et al.* (1997) used the pointwise approach in order to develop minimax theory on the Besov balls. All results in Lepski *et al.* (1997) were obtained for the one-dimensional case $d = 1$ and the selection rule proposed there, being a modification of Lepski's method, cannot be directly extended to dimensions greater than one. Generalization to an arbitrary dimension was proposed in Kerkyacharian *et al.* (2001). This allowed minimax results to be developed for the anisotropic Besov-type functional classes. The class studied in this section can be viewed as a particular case of the anisotropic one and in Theorem 4, we reproduce the results from Lepski *et al.* (1997).

### 4.5. Mixture of problems

Consider the family of estimators

$$\mathcal{F}(\mathcal{K}_{\mathcal{H},\mathcal{E}}) = \left\{\hat{F}_{h,E} : \hat{F}_{h,E}(x) = \int G^*_{h,E}(t,x) Y(\mathrm{d}t), h \in \mathcal{H}, E \in \mathcal{E}\right\}$$

and let $\hat{F}$ be the estimator derived from the collection $\mathcal{F}(\mathcal{K}_{\mathcal{H},\mathcal{E}})$ in accordance with our general selection rule, where the majorant $Q$ is given by

$$Q(\sigma) = \mathcal{C}_2 \sigma \sqrt{1 + \ln(1/\varepsilon)}.$$

Here, $\mathcal{C}_2 = \mathcal{C}_2(d,r,g)$ is a numerical constant.

**Theorem 5.** *Choose $h_{\min} = \varepsilon^2$, $h_{\max} = 1/2$ and suppose that $G$ satisfies assumption (32). Then, for all $\varepsilon > 0$ small enough,*

1. *under the conditions of Theorem 2 the estimator $\hat{F}$ is minimax, that is, it satisfies (39);*
2. *under the conditions of Theorem 3, we have*

$$\sup_{F \in \mathbb{H}_d(\boldsymbol{\alpha},L)} \mathcal{R}_r[\hat{F}; F] \leq C L^{1/(2\gamma+1)} \left(\varepsilon \sqrt{\ln \frac{1}{\varepsilon}}\right)^{2\gamma/(2\gamma+1)},$$

*where $C$ depends only on $g$, $d$, $r$ and $\gamma$;*
3. *under the conditions of Theorem 4, we have*

$$\mathcal{R}_{\mathbb{L}_r}(\hat{F}) \leq C \begin{cases} [\varepsilon \sqrt{\ln 1/\varepsilon}]^{s/(s+d/2)}, & \text{if } sp > \dfrac{d(r-p)}{2}, \\ [\varepsilon \sqrt{\ln 1/\varepsilon}]^{s/(s+d/2)} [\ln 1/\varepsilon]^{1/r}, & \text{if } sp = \dfrac{d(r-p)}{2}, \\ [\varepsilon \sqrt{\ln 1/\varepsilon}]^{s-d(1/p-1/r)/(s-d(1/p-1/2))}, & \text{if } sp < \dfrac{d(r-p)}{2}, \end{cases}$$

*where $C$ depends only on $g$, $d$, $r$ and $\gamma$.*

The proof of the theorem is along the same lines as the proofs of Theorems 2–4 and is hence omitted.

***Remark 12.***

1. Comparing the results from Theorems 3 and 5, we conclude that the rate provided by the estimator $\hat{F}$ differs from the minimax rate of convergence on $\mathbb{F}_{AH}$ by a $[\ln(1/\varepsilon)/\ln\ln(1/\varepsilon)]^{2\gamma/(2\gamma+1)}$-factor.
2. Comparing the results from Theorems 4 and 5, we conclude that the estimator $\hat{F}$ is minimax adaptive up to a $\sqrt{\ln(1/\varepsilon)}$-factor for all values of parameters $s$ and $p$. Moreover, $\hat{F}$ is a minimax adaptive estimator on the isotropic Besov balls of functions for all $s$ and $p$ such that $sp < \frac{d(r-p)}{2}$. A wavelet thresholding estimator which is nearly minimax adaptive over a scale of one-dimensional Besov balls was developed in Donoho *et al.* (1995).

## 5. Proof of Theorem 1

In the proof below, $c, c_1, c_2, \ldots$ denote constants depending only on $r$, $c_e$ and $C_e$; they can be different on different occasions.

$0^0$. We begin the proof by showing that under the premise of the theorem, the selection rule (22)–(24) is well defined, that is, there exists a measurable choice of $\hat{\mu} \in \Theta$ such that (23) is fulfilled.

It follows from Lemma 1 and assumptions (K2) and (K0) that there exists a separable modification of the Gaussian random process $\{\xi_{\mu,\nu}(x) - \xi_\nu(x), (\mu,\nu) \in \Theta \times \Theta\}$ that with probability one belongs to the $2m$-dimensional isotropic Hölder space with regularity index $0 < \tau < \gamma$ (see Lifshits (1995), Section 15). In addition, if (K2) holds and $F$ is uniformly bounded, then the integral $\int K_\nu(y,x) B_\mu(y)\,dy$, considered as a function of $(\mu,\nu)$, belongs to the $2m$-dimensional Hölder space with regularity index $\gamma$. Then, by (14) and (13), we obtain that $|\hat{F}_{\mu,\nu}(x) - \hat{F}_\nu(x)|$ is continuous in $(\mu,\nu)$. It also follows from (26) and (27) that $\sigma_\nu(x)$ (and $\tilde{\sigma}_\nu(x)$) are continuous functions of $\nu \in \Theta$. Hence, $Q(\tilde{\sigma}_\nu)$ is also continuous in $\nu$; thus, the random function under the supremum on the RHS of (22) is continuous in $(\mu,\nu)$. $\hat{R}_\mu$ is then a random variable for every $\mu \in \Theta$.

We now describe the construction of the measurable choice $\hat{\mu} \in \Theta$ satisfying (23). Let

$$\hat{R}_{\mu,\nu} = |\hat{F}_{\mu,\nu} - \hat{F}_\nu| - \tfrac{1}{2}\varepsilon Q(\tilde{\sigma}_\nu), \qquad \mu,\nu \in \Theta.$$

For any $\delta > 0$, there exists a simple function, say $\tilde{R}_{\mu,\nu}$, on $\Theta \times \Theta$ such that $|\hat{R}_{\mu,\nu} - \tilde{R}_{\mu,\nu}| \leq \delta$ for all $\mu,\nu \in \Theta$. Then, clearly,

$$|\hat{R}_\mu - \tilde{R}_\mu| \leq \delta \qquad \forall \mu \in \Theta, \tag{43}$$

where we have defined $\tilde{R}_\mu := \sup_{\nu: \tilde{\sigma}_\nu \geq \tilde{\sigma}_\mu} \tilde{R}_{\mu,\nu}$. We now observe that $\tilde{R}_\mu$ is a simple function of $\mu \in \Theta$ and define $\hat{\mu} = \arg\inf_{\mu \in \Theta}\{\tilde{R}_\mu + \varepsilon Q(\tilde{\sigma}_\mu)\}$. Since the function $\tilde{R}_\mu$ assumes a finite number of values and $Q(\tilde{\sigma}_\mu)$ is continuous in $\mu$, $\hat{\mu}$ is measurable and belongs to $\Theta$. (43) then implies (23) if $\delta$ is chosen to be $\tfrac{1}{4}\varepsilon Q(\sigma_{\min})$.

$1^0$. We write

$$\{\mathbb{E}_F |\hat{F} - F|^r\}^{1/r} \leq \{\mathbb{E}_F |\hat{F} - F|^r \mathbf{1}(\tilde{\sigma}_{\hat{\mu}} \leq \tilde{\sigma}_{\mu^*})\}^{1/r} + \{\mathbb{E}_F |\hat{F} - F|^r \mathbf{1}(\tilde{\sigma}_{\hat{\mu}} > \tilde{\sigma}_{\mu^*})\}^{1/r} \tag{44}$$

and our current goal is to bound the two terms on the right-hand side.

$2^0$. By the triangle inequality,

$$|\hat{F}_{\hat{\mu}} - F|\mathbf{1}(\tilde{\sigma}_{\hat{\mu}} \leq \tilde{\sigma}_{\mu^*}) \leq [|\hat{F}_{\hat{\mu},\mu^*} - \hat{F}_{\hat{\mu}}| + |\hat{F}_{\hat{\mu},\mu^*} - \hat{F}_{\mu^*}| + |\hat{F}_{\mu^*} - F|]\mathbf{1}(\tilde{\sigma}_{\hat{\mu}} \leq \tilde{\sigma}_{\mu^*})$$
$$=: [J_1 + J_2 + J_3]\mathbf{1}(\tilde{\sigma}_{\hat{\mu}} \leq \tilde{\sigma}_{\mu^*}).$$

By the definitions of $\hat{F}_{\mu,\nu}$, $\hat{F}_\mu$ and $\tilde{B}_\mu$,

$$J_1\,\mathbf{1}(\tilde{\sigma}_{\hat{\mu}} \leq \tilde{\sigma}_{\mu^*}) \leq [|B_{\hat{\mu},\mu^*} - B_{\hat{\mu}}| + \varepsilon|\xi_{\hat{\mu},\mu^*} - \xi_{\hat{\mu}}|]\mathbf{1}(\tilde{\sigma}_{\hat{\mu}} \leq \tilde{\sigma}_{\mu^*})$$
$$\leq \sup_{\nu \in \Theta}|B_{\nu,\mu^*} - B_\nu| + \varepsilon \sup_{\nu:\tilde{\sigma}_\nu \leq \tilde{\sigma}_{\mu^*}}|\xi_{\nu,\mu^*} - \xi_\nu|$$
$$\leq \tilde{B}_{\mu^*} + \varepsilon \sup_{\nu:\tilde{\sigma}_\nu \leq \tilde{\sigma}_{\mu^*}}|\xi_{\nu,\mu^*} - \xi_\nu|.$$

Therefore, in view of Lemma A.1 in the Appendix,

$$\{\mathbb{E}_F J_1^r \mathbf{1}(\tilde{\sigma}_{\hat{\mu}} \leq \tilde{\sigma}_{\mu^*})\}^{1/r} \leq \tilde{B}_{\mu^*} + \varepsilon\bigg\{\mathbb{E}\sup_{\nu:\tilde{\sigma}_\nu \leq \tilde{\sigma}_{\mu^*}}|\xi_{\nu,\mu^*} - \xi_\nu|^r\bigg\}^{1/r}$$
$$\leq \tilde{B}_{\mu^*} + C_r \varepsilon\{e(\tilde{\sigma}_{\mu^*}) + 2\tilde{\sigma}_{\mu^*}\} \qquad (45)$$
$$\leq c\varepsilon Q(\tilde{\sigma}_{\mu^*}),$$

where we have used the definitions of $\mu^*$ and $Q(\cdot)$.

Furthermore,

$$J_2\mathbf{1}(\tilde{\sigma}_{\hat{\mu}} \leq \tilde{\sigma}_{\mu^*}) \leq \{|\hat{F}_{\hat{\mu},\mu^*} - \hat{F}_{\mu^*}| - \tfrac{1}{2}\varepsilon Q(\tilde{\sigma}_{\mu^*})\}\mathbf{1}(\tilde{\sigma}_{\hat{\mu}} \leq \tilde{\sigma}_{\mu^*}) + \tfrac{1}{2}\varepsilon Q(\tilde{\sigma}_{\mu^*})$$
$$\leq \hat{R}_{\hat{\mu}} + \tfrac{1}{2}\varepsilon Q(\tilde{\sigma}_{\mu^*})$$
$$\leq \hat{R}_{\mu^*} + \tfrac{3}{2}\varepsilon Q(\tilde{\sigma}_{\mu^*}) + \delta,$$

where the second inequality follows from (22) and the third is a consequence of (23). Hence,

$$\{\mathbb{E}_F J_2^r \,\mathbf{1}(\tilde{\sigma}_{\hat{\mu}} \leq \tilde{\sigma}_{\mu^*})\}^{1/r} \leq \{\mathbb{E}_F \hat{R}_{\mu^*}^r \mathbf{1}(\hat{R}_{\mu^*} > 0)\}^{1/r} + \tfrac{3}{2}\varepsilon Q(\tilde{\sigma}_{\mu^*}) + \delta.$$

Because

$$\hat{R}_{\mu^*} = \sup_{\nu:\tilde{\sigma}_\nu \geq \tilde{\sigma}_{\mu^*}}[|\hat{F}_{\mu^*,\nu} - \hat{F}_\nu| - \tfrac{1}{2}\varepsilon Q(\tilde{\sigma}_\nu)]$$
$$\leq \tilde{B}_{\mu^*} + \varepsilon\bigg[\sup_{\nu:\tilde{\sigma}_\nu \geq \tilde{\sigma}_{\mu^*}}|\xi_{\mu^*,\nu} - \xi_\nu| - \tfrac{1}{2}Q(\tilde{\sigma}_\nu)\bigg]_+, \qquad ([\cdot]_+ = \max\{\cdot,0\}), \qquad (46)$$

we obtain

$$\{\mathbb{E}_F J_2^r \,\mathbf{1}(\tilde{\sigma}_{\hat{\mu}} \leq \tilde{\sigma}_{\mu^*})\}^{1/r}$$

$$\leq \tilde{B}_{\mu^*} + \varepsilon \left\{ \mathbb{E} \left[ \sup_{\nu : \tilde{\sigma}_\nu \geq \tilde{\sigma}_{\mu^*}} |\xi_{\mu^*,\nu} - \xi_\nu| - \tfrac{1}{2} Q(\tilde{\sigma}_\nu) \right]_+^r \right\}^{1/r} + \tfrac{3}{2} \varepsilon Q(\tilde{\sigma}_{\mu^*}) + \delta$$
$$\leq \tilde{B}_{\mu^*} + \tfrac{3}{2} \varepsilon Q(\tilde{\sigma}_{\mu^*}) + \delta + c\varepsilon \sigma_{\min} \qquad (47)$$
$$\leq c\varepsilon Q(\tilde{\sigma}_{\mu^*}),$$

where the second inequality follows from Lemma A.2 (see Appendix) and the last inequality follows from the definitions of $\mu^*$ and $Q(\cdot)$.

The bound on $\{\mathbb{E}_F J_3^r \mathbf{1}(\tilde{\sigma}_{\hat{\mu}} \leq \tilde{\sigma}_{\mu^*})\}^{1/r}$ is immediate:

$$\{\mathbb{E}_F J_3^r \mathbf{1}(\tilde{\sigma}_{\hat{\mu}} \leq \tilde{\sigma}_{\mu^*})\}^{1/r} \leq \tilde{B}_{\mu^*} + \varepsilon [\mathbb{E}|\xi_{\mu^*}|^r]^{1/r} \leq \tilde{B}_{\mu^*} + c\varepsilon \sigma_{\mu^*} \leq c\varepsilon Q(\tilde{\sigma}_{\mu^*}). \qquad (48)$$

Combining (45), (47) and (48), we obtain that there exists a constant $c$ depending only on $r$ such that

$$\{\mathbb{E}_F |\hat{F}_{\hat{\mu}} - F|^r \mathbf{1}(\tilde{\sigma}_{\hat{\mu}} \leq \tilde{\sigma}_{\mu^*})\}^{1/r} \leq c\varepsilon Q(\tilde{\sigma}_{\mu^*}). \qquad (49)$$

$3^0$. To bound the second term on the right-hand side of (44), we proceed as follows. Define the events $A_k = \{2^{k-1} \tilde{\sigma}_{\mu^*} \leq \tilde{\sigma}_{\hat{\mu}} < 2^k \tilde{\sigma}_{\mu^*}\}$, $k = 1, 2, \ldots$, and let $\mu_k \in \Theta_F(\mathcal{K}_\Theta)$ be such that the corresponding estimators $\hat{F}_{\mu_k} \in \mathcal{F}(\mathcal{K}_\Theta)$ have the following properties:

(i) $\tilde{\sigma}_{\mu_k}^2 = \mathrm{var}\{\hat{F}_{\mu_k}\} = 2^k \tilde{\sigma}_{\mu^*}$;
(ii) $\tilde{B}_{\mu_k} \leq \tfrac{1}{4} \varepsilon Q(\tilde{\sigma}_{\mu_k})$.

The existence of estimators $\hat{F}_{\mu_k}$ satisfying (i) and (ii) is guaranteed by the fact that $\Theta_F(\mathcal{K}_\Theta)$ is non-empty and $\mu^* \in \Theta_F(\mathcal{K}_\Theta)$. We can then write

$$|\hat{F}_{\hat{\mu}} - F| \mathbf{1}(\tilde{\sigma}_{\hat{\mu}} \geq \tilde{\sigma}_{\mu^*}) \leq \sum_{k=1}^\infty [|\hat{F}_{\hat{\mu},\mu_k} - \hat{F}_{\hat{\mu}}| + |\hat{F}_{\hat{\mu},\mu_k} - \hat{F}_{\mu_k}| + |\hat{F}_{\mu_k} - F|] \mathbf{1}(A_k)$$
$$=: \sum_{k=1}^\infty [I_{1,k} + I_{2,k} + I_{3,k}] \mathbf{1}(A_k). \qquad (50)$$

We have

$$I_{1,k} \mathbf{1}(A_k) \leq [\tilde{B}_{\mu_k} + \varepsilon |\xi_{\hat{\mu},\mu_k} - \xi_{\hat{\mu}}|] \mathbf{1}(A_k)$$
$$\leq \left[ \tfrac{1}{2} \varepsilon Q(\tilde{\sigma}_{\mu_k}) + \varepsilon \sup_{\nu : \tilde{\sigma}_\nu \leq \tilde{\sigma}_{\mu_k}} |\xi_{\nu,\mu_k} - \xi_\nu| \right] \mathbf{1}(A_k),$$

where the second inequality follows from the definition of $\mu_k$. Hence, by the Cauchy–Schwarz inequality and Lemma A.1,

$$\{\mathbb{E}_F I_{1,k}^r \mathbf{1}(A_k)\}^{1/r}$$

$$\leq \tfrac{1}{2}\varepsilon Q(\tilde{\sigma}_{\mu_k})\mathbb{P}_F^{1/r}(A_k) + \varepsilon \bigg\{\mathbb{E}_F \sup_{\nu:\tilde{\sigma}_\nu \leq \tilde{\sigma}_{\mu_k}} |\xi_{\nu,\mu_k} - \xi_\nu|^{2r}\bigg\}^{1/2r}[\mathbb{P}_F(A_k)]^{1/2r}$$

$$\leq \tfrac{1}{2}\varepsilon Q(\tilde{\sigma}_{\mu_k})\mathbb{P}_F^{1/r}(A_k) + c\varepsilon\{e(\tilde{\sigma}_{\mu_k}) + \tilde{\sigma}_{\mu_k}\}[\mathbb{P}_F(A_k)]^{1/2r} \tag{51}$$

$$\leq c\varepsilon Q(\tilde{\sigma}_{\mu_k})[\mathbb{P}_F(A_k)]^{1/2r}.$$

Furthermore,

$$\begin{aligned}
I_{2,k}\mathbf{1}(A_k) &\leq [|\hat{F}_{\hat{\mu},\mu_k} - \hat{F}_{\mu_k}| - \tfrac{1}{2}\varepsilon Q(\tilde{\sigma}_{\mu_k})]\mathbf{1}(A_k) + \tfrac{1}{2}\varepsilon Q(\tilde{\sigma}_{\mu_k})\mathbf{1}(A_k)\\
&\leq [\hat{R}_{\hat{\mu}} + \tfrac{1}{2}\varepsilon Q(\tilde{\sigma}_{\mu_k})]\mathbf{1}(A_k)\\
&\leq [\hat{R}_{\mu^*} + \tfrac{3}{2}\varepsilon Q(\tilde{\sigma}_{\mu_k}) + \delta]\mathbf{1}(A_k)\\
&\leq [\hat{R}_{\mu^*}\mathbf{1}(\hat{R}_{\mu^*} > 0) + \tfrac{3}{2}\varepsilon Q(\tilde{\sigma}_{\mu_k}) + \delta]\mathbf{1}(A_k),
\end{aligned}$$

where the second inequality follows from the definition of $\hat{R}_\mu$ and the third from the definition of $\hat{\mu}$ and the monotonicity of $Q(\cdot)$. Arguing as in (46) and (47), and using the Cauchy–Schwarz inequality, we obtain

$$\begin{aligned}
\{\mathbb{E}_F I_{2,k}^r \mathbf{1}(A_k)\}^{1/r} &\\
&\leq \{\mathbb{E}_F \hat{R}_{\mu^*}^{2r}\mathbf{1}(\hat{R}_{\mu^*} > 0)\}^{1/2r}[\mathbb{P}_F(A_k)]^{1/2r} + (\tfrac{3}{2}\varepsilon Q(\tilde{\sigma}_{\mu_k}) + \delta)[\mathbb{P}_F(A_k)]^{1/r}\\
&\leq (\tilde{B}_{\mu^*} + c\varepsilon\sigma_{\min})[\mathbb{P}_F(A_k)]^{1/2r} + (\tfrac{3}{2}\varepsilon Q(\tilde{\sigma}_{\mu_k}) + \delta)[\mathbb{P}_F(A_k)]^{1/r}\\
&\leq c\varepsilon Q(\tilde{\sigma}_{\mu_k})\}[\mathbb{P}_F(A_k)]^{1/2r}.
\end{aligned} \tag{52}$$

Finally,

$$\begin{aligned}
\{\mathbb{E}_F I_{3,k}^r \mathbf{1}(A_k)\}^{1/r} &\leq \tilde{B}_{\mu_k}[\mathbb{P}_F(A_k)]^{1/r} + \varepsilon\{\mathbb{E}|\xi_{\mu_k}|^{2r}\}^{1/2r}[\mathbb{P}_F(A_k)]^{1/2r}\\
&\leq \tfrac{1}{2}\varepsilon Q(\tilde{\sigma}_{\mu_k})[\mathbb{P}_F(A_k)]^{1/r} + c\varepsilon\tilde{\sigma}_{\mu_k}[\mathbb{P}_F(A_k)]^{1/2r}\\
&\leq c\varepsilon Q(\tilde{\sigma}_{\mu_k})[\mathbb{P}_F(A_k)]^{1/2r},
\end{aligned} \tag{53}$$

where we have used the definition of $\mu_k$. Combining (51), (52) and (53), we obtain

$$\{\mathbb{E}_F I_{1,k}^r\mathbf{1}(A_k)\}^{1/r} + \{\mathbb{E}_F I_{2,k}^r\mathbf{1}(A_k)\}^{1/r} + \{\mathbb{E}_F I_{3,k}^r\mathbf{1}(A_k)\}^{1/r} \leq c\varepsilon Q(\tilde{\sigma}_{\mu_k})\}[\mathbb{P}_F(A_k)]^{1/2r}. \tag{54}$$

In order to complete the proof, we need to bound $\mathbb{P}_F(A_k)$ from above.

$4^0$. Note that for any integer $1 < m < k$, by definition of $\hat{\mu}$, we have $A_k \subseteq \{\tilde{\sigma}_{\hat{\mu}} > \tilde{\sigma}_{\mu_{k-m}}\}$. Hence,

$$\begin{aligned}
A_k &\subseteq \{\tilde{\sigma}_{\hat{\mu}} > \tilde{\sigma}_{\mu_{k-m}}\}\\
&\subseteq \{\hat{R}_{\hat{\mu}} + \varepsilon Q(\tilde{\sigma}_{\hat{\mu}}) < \hat{R}_{\mu_{k-m}} + \varepsilon Q(\tilde{\sigma}_{\mu_{k-m}}) + \delta\}
\end{aligned} \tag{55}$$

$$\subseteq \{\hat{R}_{\mu_{k-m}} + \varepsilon Q(\tilde{\sigma}_{\mu_{k-m}}) > \tfrac{1}{2}\varepsilon Q(\tilde{\sigma}_{\hat{\mu}}) - \delta\}$$
$$\subseteq \{\hat{R}_{\mu_{k-m}} + \varepsilon Q(\tilde{\sigma}_{\mu_{k-m}}) > \tfrac{1}{2}\varepsilon Q(\tilde{\sigma}_{\mu_{k-1}}) - \delta\},$$

where the second inclusion is by (25) and the third is by the monotonicity of $Q(\cdot)$. Furthermore, using assumption (E), we have

$$\frac{1}{2}Q(\tilde{\sigma}_{\mu_{k-1}}) - Q(\tilde{\sigma}_{\mu_{k-m}})$$

$$= \frac{1}{2}\varkappa_0 e(\tilde{\sigma}_{\mu_{k-1}}) + \frac{1}{2}\tilde{\sigma}_{\mu_{k-1}}\sqrt{1 + \varkappa_1 \ln \frac{\tilde{\sigma}_{\mu_{k-1}}}{\sigma_{\min}}} - \varkappa_0 e(\tilde{\sigma}_{\mu_{k-m}}) - \tilde{\sigma}_{\mu_{k-m}}\sqrt{1 + \varkappa_1 \ln \frac{\tilde{\sigma}_{\mu_{k-m}}}{\sigma_{\min}}}$$

$$\geq \left[\frac{1}{2}c_e^{m-1} - 1\right]\varkappa_0 e(\tilde{\sigma}_{\mu_{k-m}}) + \tilde{\sigma}_{\mu_{k-2}}\sqrt{1 + \varkappa_1 \ln \frac{\tilde{\sigma}_{\mu_{k-2}}}{\sigma_{\min}}} - \tilde{\sigma}_{\mu_{k-m}}\sqrt{1 + \ln \frac{\tilde{\sigma}_{\mu_{k-m}}}{\sigma_{\min}}}$$

$$\geq \frac{1}{2}\varkappa_0 e(\tilde{\sigma}_{\mu_{k-m}}) + \frac{1}{2}\tilde{\sigma}_{\mu_{k-m}}\sqrt{1 + \varkappa_1 \ln \frac{\tilde{\sigma}_{\mu_{k-m}}}{\sigma_{\min}}}$$

$$= \frac{1}{2}Q(\tilde{\sigma}_{\mu_{k-m}}),$$

provided that $m \geq 3 \vee [1 + (\ln 3/\ln c_e)]$. Choosing

$$m = m_0 := \lceil 1 + (\ln 3/\ln c_e) \rceil \vee 3, \tag{56}$$

we obtain that

$$\mathbb{P}_F(A_k) \leq \mathbb{P}_F\left\{\hat{R}_{\mu_k-m_0} > \frac{1}{2}\varepsilon Q(\tilde{\sigma}_{\mu_k-m_0}) - \delta\right\}$$

$$\leq \mathbb{P}_F\bigg\{\tilde{B}_{\mu_k-m_0}$$

$$+ \varepsilon \sup_{\nu:\tilde{\sigma}_\nu \geq \tilde{\sigma}_{\mu_k-m_0}}\left[|\xi_{\mu_k-m_0,\nu} - \xi_\nu| - \frac{1}{2}Q(\tilde{\sigma}_{\mu_k-m_0})\right] > \frac{1}{4}\varepsilon Q(\tilde{\sigma}_{\mu_k-m_0})\bigg\} \tag{57}$$

$$\leq \mathbb{P}_F\bigg\{\sup_{\nu:\tilde{\sigma}_\nu \geq \tilde{\sigma}_{\mu_k-m_0}}\left[|\xi_{\mu_k-m_0,\nu} - \xi_\nu| - \frac{1}{2}\varepsilon Q(\tilde{\sigma}_{\mu_k-m_0})\right] > 0\bigg\}$$

$$\leq 4\left(\frac{\sigma_{\min}}{\tilde{\sigma}_{\mu_k-m_0}}\right)^{\varkappa_1/64},$$

where the first inequality is by (55), the second is by the bound on $\hat{R}_\mu$ (see (46)), the definition of $\delta$ and the monotonicity of $Q(\cdot)$, the third is in view of the definition of the $\mu_k$'s and the fourth inequality follows from Lemma A.2.

$5^0$. Now using (54) and (57), we bound $\{\mathbb{E}_F|\hat{F}_{\hat{\mu}} - F|^r \mathbf{1}(\tilde{\sigma}_{\hat{\mu}} \geq \tilde{\sigma}_{\mu^*})\}^{1/r}$; see (50).

Let $m_0$ be given by (56) and, for the sake of brevity, set $\gamma = \varkappa_1/64$. Then,

$$\sum_{k=m_0+1}^{\infty} [\mathbb{P}_F(A_k)]^{1/2r} \leq c\left(\frac{\sigma_{\min}}{\tilde{\sigma}_{\mu^*}}\right)^{\gamma/(2r)} 2^{m_0\gamma/(2r)} \sum_{k=m_0+1}^{\infty} 2^{-k\gamma/(2r)} \leq c\left(\frac{\sigma_{\min}}{\tilde{\sigma}_{\mu^*}}\right)^{\gamma/(2r)} 2^{m_0\gamma/(2r)}.$$

Moreover, using assumption (E), we obtain

$$\sum_{k=m_0+1}^{\infty} Q(\tilde{\sigma}_{\mu_k})[\mathbb{P}_F(A_k)]^{1/2r}$$

$$\leq c \sum_{k=m_0+1}^{\infty} \left[\varkappa_0 e(\tilde{\sigma}_{\mu_k}) + \tilde{\sigma}_{\mu_k}\sqrt{1 + \varkappa_1 \ln \frac{\tilde{\sigma}_{\mu_k}}{\sigma_{\min}}}\right]\left(\frac{\sigma_{\min}}{\tilde{\sigma}_{\mu_{k-m_0}}}\right)^{\gamma/(2r)}$$

$$\leq c_1 \varkappa_0 2^{m_0\gamma/(2r)}\left(\frac{\sigma_{\min}}{\tilde{\sigma}_{\mu^*}}\right)^{\gamma/(2r)} e(\tilde{\sigma}_{\mu^*}) \sum_{k=m_0+1}^{\infty} C_e^k 2^{-k\gamma/(2r)}$$

$$+ c_2 2^{m_0\gamma/(2r)} \tilde{\sigma}_{\mu^*}\left(\frac{\sigma_{\min}}{\tilde{\sigma}_{\mu^*}}\right)^{\gamma/(2r)} \sum_{k=m_0+1}^{\infty} 2^{k-k\gamma/(2r)}\sqrt{1 + \varkappa_1 \ln \frac{2^k \tilde{\sigma}_{\mu^*}}{\sigma_{\min}}} \quad (58)$$

$$\leq c 2^{m_0\gamma/(2r)}\left(\frac{\sigma_{\min}}{\tilde{\sigma}_{\mu^*}}\right)^{\gamma/(2r)}\left(\varkappa_0 e(\tilde{\sigma}_{\mu^*}) + \tilde{\sigma}_{\mu^*}\sqrt{1 + \varkappa_1 \ln \frac{\tilde{\sigma}_{\mu^*}}{\sigma_{\min}}}\right),$$

$$\leq c 2^{m_0\gamma/(2r)}\left(\frac{\sigma_{\min}}{\tilde{\sigma}_{\mu^*}}\right)^{\gamma/(2r)} Q(\tilde{\sigma}_{\mu^*})$$

because, by our choice of $\varkappa_1$, $\gamma = \varkappa_1/64 \geq 2r(\ln C_e/\ln 2)$, which implies that the sums on the right-hand side are finite. In addition,

$$\sum_{k=1}^{m_0} Q(\tilde{\sigma}_{\mu_k}) = \sum_{k=1}^{m_0}\left[\varkappa_0 e(\tilde{\sigma}_{\mu_k}) + \tilde{\sigma}_{\mu_k}\sqrt{1 + \varkappa_1 \ln \frac{\tilde{\sigma}_{\mu_k}}{\sigma_{\min}}}\right]$$

$$\leq \varkappa_0 e(\tilde{\sigma}_{\mu^*}) \sum_{k=1}^{m_0} C_e^k + \tilde{\sigma}_{\mu^*}\sqrt{1 + \varkappa_1 \ln \frac{\tilde{\sigma}_{\mu^*}}{\sigma_{\min}}} \sum_{k=1}^{m_0} 2^k + 2\tilde{\sigma}_{\mu^*}\sqrt{\varkappa_1 \ln 2}\sum_{k=1}^{m_0} 2^k \sqrt{k} \quad (59)$$

$$\leq cQ(\tilde{\sigma}_{\mu^*});$$

here, we have used assumption (E).

Therefore, combining (54), (50), (58) and (59), we finally obtain

$$\{\mathbb{E}_F|\hat{F}_{\hat{\mu}} - F|^r \mathbf{1}(\tilde{\sigma}_{\hat{\mu}} \geq \tilde{\sigma}_{\mu^*})\}^{1/r}$$

$$\leq c_1\varepsilon \sum_{k=1}^{m_0} Q(\tilde{\sigma}_{\mu_k}) + c_2\varepsilon \sum_{k=m_0+1}^{\infty} Q(\tilde{\sigma}_{\mu_k})[\mathbb{P}_F(A_k)]^{1/2r} \leq c\varepsilon Q(\tilde{\sigma}_{\mu^*}).$$

This inequality and (49) lead to the statement of the theorem.

# 6. Proofs of Theorems 2, 3 and 4

The proofs of Theorems 2, 3 and 4 use upper bounds on the function $e_{\mathcal{K}_\Theta}(\cdot)$ defined in (19). Therefore, we begin this section with two lemmas establishing such bounds. We then present the proofs of Theorems 2, 3 and 4.

## 6.1. Bounds on function $e_{\mathcal{K}_\Theta}(\cdot)$

For fixed $x \in \mathcal{D}_0$, consider the random process $\{\eta_{\mu,\nu}(x), \mu, \nu \in \Theta\}$ given by

$$\eta_{\mu,\nu}(x) = \xi_{\mu,\nu}(x) - \xi_\nu(x) = \int [K_{\mu,\nu}(t,x) - K_\nu(t,x)] W(dt), \qquad \mu, \nu \in \Theta.$$

For $\lambda, \lambda' \in \Theta$, define

$$\overline{\rho}(\lambda, \lambda') = \|\widetilde{K}_\lambda - \widetilde{K}_{\lambda'}\|_{2,\infty}, \qquad \widetilde{K}_\lambda(\cdot, x) = K_\lambda(\cdot, x) / \|K_\lambda(\cdot, x)\|_2;$$
$$\underline{\rho}(\lambda, \lambda') = \sup_x |1 - \|K_\lambda(\cdot, x)\|_2 / \|K_{\lambda'}(\cdot, x)\|_2|.$$

The next lemma establishes an upper bound on the intrinsic semi-metric of the process $\{\eta_{\mu,\nu}(x), \mu, \nu \in \Theta\}$.

**Lemma 1.** *Let*

$$\rho[(\mu,\nu),(\mu',\nu')] := \sqrt{\mathbb{E}|\eta_{\mu,\nu}(x) - \eta_{\mu',\nu'}(x)|^2}.$$

(i) *Then, for all $\mu, \nu, \mu', \nu' \in \Theta$, we have*

$$\rho[(\mu,\nu),(\mu',\nu')] \le 2\tilde{\sigma}_\nu(x)[\overline{\rho}(\nu,\nu') + \underline{\rho}(\nu,\nu')] + \tilde{\sigma}_\mu(x)[\overline{\rho}(\mu,\mu') + \underline{\rho}(\mu,\mu')]. \qquad (60)$$

(ii) *In addition, suppose that $\Theta = \bigotimes_{j=1}^l \Theta_j$ and $\theta = (\theta_1, \ldots, \theta_l)$, where $\theta_j \in \Theta_j$, $j = 1, \ldots, l$. Given $\vec{\theta}_j \in \bigotimes_{i=1, i \ne j}^l \Theta_i$ and $\lambda_j \in \Theta_j$, write $K_{\vec{\theta}_j, \lambda_j} = K_{\theta_1, \ldots, \theta_{j-1}, \lambda_j, \theta_{j+1}, \ldots, \theta_l}$. Then,*

$$\overline{\rho}(\lambda, \lambda') \le \sum_{\substack{j=1 \\ \lambda_j \ne \lambda'_j}}^l \sup_{\vec{\theta}_j} \|\widetilde{K}_{\vec{\theta}_j, \lambda_j} - \widetilde{K}_{\vec{\theta}_j, \lambda'_j}\|_{2,\infty} \qquad \forall \lambda, \lambda' \in \Theta. \qquad (61)$$

**Proof.** (i) We have

$$\rho[(\mu,\nu),(\mu',\nu')] = \|(K_{\nu,\mu}(\cdot,x) - K_\nu(\cdot,x)) - (K_{\nu',\mu'}(\cdot,x) - K_{\nu'}(\cdot,x))\|_2$$
$$\le \|K_{\nu,\mu}(\cdot,x) - K_{\nu',\mu'}(\cdot,x)\|_2 + \|K_\nu(\cdot,x) - K_{\nu'}(\cdot,x)\|_2$$
$$\le \|K_{\nu,\mu}(\cdot,x) - K_{\nu',\mu}(\cdot,x)\|_2 + \|K_{\nu',\mu}(\cdot,x) - K_{\nu',\mu'}(\cdot,x)\|_2$$

$$+ \|K_\nu(\cdot, x) - K_{\nu'}(\cdot, x)\|_2$$
$$= \|K_{\nu,\mu}(\cdot, x) - K_{\nu',\mu}(\cdot, x)\|_2 + \|K_{\mu,\nu'}(\cdot, x) - K_{\mu',\nu'}(\cdot, x)\|_2$$
$$+ \|K_\nu(\cdot, x) - K_{\nu'}(\cdot, x)\|_2,$$

where the last line follows from assumption (K1).

Thus, to prove (60), it suffices to show that for all $\lambda, \lambda' \in \Theta$,

$$\sup_{\theta \in \Theta} \|K_{\lambda,\theta}(\cdot, x) - K_{\lambda',\theta}(\cdot, x)\|_2 \leq \tilde{\sigma}_\lambda(x)(\overline{\rho}(\lambda, \lambda') + \underline{\rho}(\lambda, \lambda')). \tag{62}$$

Let us prove (62). Indeed, using the Minkowski inequality, we get, for all $\theta \in \Theta$,

$$\|K_{\lambda,\theta}(\cdot, x) - K_{\lambda',\theta}(\cdot, x)\|_2 = \sqrt{\int \left(\int K_\theta(y, x)[K_\lambda(t, y) - K_{\lambda'}(t, y)] \,\mathrm{d}y\right)^2 \mathrm{d}t}$$
$$\leq \int |K_\theta(y, x)| \|K_\lambda(\cdot, y) - K_{\lambda'}(\cdot, y)\|_2 \,\mathrm{d}y.$$

Moreover, for all $y$,

$$\|K_\lambda(\cdot, y) - K_{\lambda'}(\cdot, y)\|_2 \leq \|K_\lambda(\cdot, y)\|_2 \|\widetilde{K}_\lambda(\cdot, y) - \widetilde{K}_{\lambda'}(\cdot, y)\|_2 + |\|K_\lambda(\cdot, y)\|_2 - \|K_{\lambda'}(\cdot, y)\|_2|$$
$$\leq \|K_\lambda(\cdot, y)\|_2(\overline{\rho}(\lambda, \lambda') + \underline{\rho}(\lambda, \lambda')).$$

It remains to note that by definition,

$$\tilde{\sigma}_\nu(x) = \sup_{\theta \in \Theta} \int |K_\theta(y, x)| \|K_\nu(\cdot, y)\|_2 \,\mathrm{d}y \vee \sigma_\nu(x).$$

(ii) The statement follows immediately from the triangle inequality. □

Using general results of Lemma 1, we now establish an upper bound on the intrinsic semi-metric of the Gaussian process $\eta_{\mu,\nu}$ with index set $\Theta = \mathcal{H} \times \mathcal{E}$.

**Lemma 2.** *Let $\mathcal{K}_{\mathcal{H},\mathcal{E}}$ be the family of kernels defined in (34). Then, for all $\mu, \nu, \nu' \in \mathcal{H} \times \mathcal{E}$, we have*

$$\rho[(\mu, \nu), (\mu, \nu')] \leq 2\tilde{\sigma}_\nu(x) \left\{ M \left( \sum_{i=1}^d \left| 1 - \frac{h_i}{h'_i} \right|^2 \right)^{1/2} + 2 \left| 1 - \prod_{i=1}^d \frac{h'_i}{h_i} \right| + Mdh_{\min}^{-1}|E - E'|_2 \right\}.$$

**Proof.** Let $\nu = (h, E)$, $\nu' = (h', E') \in \mathcal{H} \times \mathcal{E}$. Our current goal is to bound $\rho(\nu, \nu')$ from above. For this purpose, we apply Lemma 1 with $\Theta = \mathcal{H} \times \mathcal{E}$.

In view of (33), (34) and Remark 7, we have

$$\|\sigma_{h,E}^{-1} G_{h,E} - \sigma_{h',E}^{-1} G_{h',E}\|_2$$

$$= \|G\|_2^{-1} \left\{ \int \left| \left( \prod_{i=1}^{d} \frac{1}{h_i^{1/2}} \right) G\left( \frac{t_1}{h_1}, \ldots, \frac{t_d}{h_d} \right) - \left( \prod_{i=1}^{d} \frac{1}{[h_i']^{1/2}} \right) G\left( \frac{t_1}{h_1'}, \ldots, \frac{t_d}{h_d'} \right) \right|^2 dt \right\}^{1/2}$$

$$= \|G\|_2^{-1} \left\{ \int \left| G(t_1, \ldots, t_d) - \left( \prod_{i=1}^{d} \sqrt{h_i/h_i'} \right) G\left( \frac{h_1}{h_1'} t_1, \ldots, \frac{h_d}{h_d'} t_d \right) \right|^2 dt \right\}^{1/2}$$

$$\leq \|G\|_2^{-1} \left\{ \int \left| G(t_1, \ldots, t_d) - G\left( \frac{h_1}{h_1'} t_1, \ldots, \frac{h_d}{h_d'} t_d \right) \right|^2 dt \right\}^{1/2} \qquad (63)$$

$$+ \|G\|_2^{-1} \left| 1 - \prod_{i=1}^{d} \sqrt{h_i/h_i'} \right| \left\{ \int G^2 \left( \frac{h_1}{h_1'} t_1, \ldots, \frac{h_d}{h_d'} t_d \right) dt \right\}^{1/2}$$

$$\leq \|G\|_2^{-1} M \left\{ \sum_{i=1}^{d} \left| 1 - \frac{h_i}{h_i'} \right|^2 \right\}^{1/2} + \prod_{i=1}^{d} \sqrt{h_i'/h_i} \left| 1 - \prod_{i=1}^{d} \sqrt{h_i/h_i'} \right|$$

$$\leq M \left\{ \sum_{i=1}^{d} \left| 1 - \frac{h_i}{h_i'} \right|^2 \right\}^{1/2} + \left| 1 - \prod_{i=1}^{d} \sqrt{h_i'/h_i} \right|;$$

here, we have taken into account (32) and the fact that $\|G\|_2 \geq 1$.

Furthermore, if $H = \mathrm{diag}(h_1, \ldots, h_d\}$, then

$$\|\sigma_{h,E}^{-1} G_{h,E} - \sigma_{h,E'}^{-1} G_{h,E'} \|_2 = \sigma_{h,E}^{-1} \|G_{h,E} - G_{h,E'}\|_2$$

$$= \|G\|_2^{-1} \prod_{i=1}^{d} h_i^{1/2} \left\{ \int |G_h(E^T t) - G_h((E')^T t)|^2 dt \right\}^{1/2}$$

$$= \|G\|_2^{-1} \left\{ \int |G(H^{-1} E^T t) - G(H^{-1}(E')^T t)|^2 dt \right\}^{1/2} \qquad (64)$$

$$\leq \|G\|_2^{-1} M \left\{ \int |H^{-1}(E - E')^T t|_2^2 dt \right\}^{1/2}$$

$$\leq M d h_{\min}^{-1} |E - E'|_2.$$

Combining (63), (64) and using (61), we obtain

$$\overline{\rho}(\nu, \nu') \leq M \left\{ \sum_{i=1}^{d} \left| 1 - \frac{h_i}{h_i'} \right|^2 \right\}^{1/2} + \left| 1 - \prod_{i=1}^{d} \sqrt{h_i'/h_i} \right| + M d h_{\min}^{-1} |E - E'|_2.$$

Observe, also, that

$$\underline{\rho}(\nu, \nu') = |1 - \sigma_{h,E}/\sigma_{h',E'}| = \left| 1 - \prod_{i=1}^{d} \sqrt{h_i'/h_i} \right|.$$

Applying Lemma 1(i), we complete the proof. □

**Lemma 3.** *There exist constants $C_1$, $C_2$ and $C_3$ depending on $G$ and $d$ such that the following statements hold:*

(i) *if $\mathcal{K}_{SI}$ is the family of kernels defined in (35) and $\sigma \geq 1$, then*

$$e_{\mathcal{K}_{SI}}(\sigma) \leq C_1 \sigma \sqrt{\ln \sigma} \qquad \forall \sigma \in \Sigma_{\Theta_{SI}};$$

(ii) *if $\mathcal{K}_{AH}$ is the family of kernels defined in (41) and $\ln(h_{\max}/h_{\min}) \geq 1$, then*

$$e_{\mathcal{K}_{AH}}(\sigma) \leq C_2 \sigma [\sqrt{\ln \ln(h_{\max}/h_{\min})} + 1] \qquad \forall \sigma \in \Sigma_{\Theta_{AH}};$$

(iii) *if $\mathcal{K}_B$ is the family of kernels defined in Section 4.4 and $\sigma \geq 1$, then*

$$e_{\mathcal{K}_B}(\sigma) \leq C_3 \sigma \sqrt{\ln(1 + \ln(\sigma/\sigma_{\min}))}.$$

**Proof.** Throughout the proof, $c, c_1, c_2, \ldots$ denote positive constants depending only on $G$ and $d$. They can be differ from appearance to appearance.

$1^0$. Let $\mathcal{K}_{SI}$ be the family of kernels defined in (35). Let $\sigma_{\min}$ and $\sigma_{\max}$ be as defined in (37) and fix $\sigma \in [\sigma_{\min}, \sigma_{\max}]$. Here, $\nu = (h, E)$ and the index set of the corresponding random process $\{\xi_{\mu,\nu} - \xi_\nu\}$ is given by $\{\nu : \tilde{\sigma}_\nu \leq \sigma\} = [h_\sigma, h_{\max}] \times \mathcal{E}$, where $h_\sigma = c_1 h_{\max}^{-d+1} \sigma^{-2}$ (see (36)).

Lemma 2 implies that the following upper bounds holds on the semi-metric $\rho_{SI}$ of this process:

$$\rho_{SI}[(\mu,\nu),(\mu,\nu')] \leq 2\sigma \left\{ M\left|1 - \frac{h_1}{h_1'}\right| + 2\left|1 - \frac{h_1'}{h_1}\right| + M d h_\sigma^{-1} |E - E'|_2 \right\},$$

for all $\nu, \nu'$ such that $\tilde{\sigma}_\nu \vee \tilde{\sigma}_{\nu'} \leq \sigma$. Note, also, that by (18),

$$\sup_{\nu : \tilde{\sigma}_\nu \leq \sigma} \text{var}(\xi_{\mu,\nu} - \xi_\nu) \leq 2 \sup_{\nu : \tilde{\sigma}_\nu \leq \sigma} \tilde{\sigma}_\nu = 2\sigma. \qquad (65)$$

The number of balls $N_1(\zeta)$ of radius $\zeta$ in semi-metric $c_2\sigma\{|1 - h_1/h_1'| + |1 - h_1'/h_1|\}$ covering the set $[h_\sigma, h_{\max}] = [c_1 h_{\max}^{-d+1} \sigma^{-2}, h_{\max}]$ admits the following upper bound:

$$N_1(\zeta) \leq \ln(c_3 \sigma^2 h_{\max}^d) \ln^{-1}(1 + c_4 \zeta \sigma^{-1}).$$

The number of balls $N_2(\zeta)$ of radius $\zeta$ in the semi-metric $c_5 \sigma |E - E'|_2$ covering $\mathcal{E}$ does not exceed $(c_6 \sigma h_\sigma^{-1} \zeta^{-1})^{d-1} = (c_7 \sigma^3 h_{\max}^{d-1} \zeta^{-1})^{d-1}$. Thus, the total number of balls covering $[h_\sigma, h_{\max}] \times \mathcal{E}$ equals $N_1(\zeta) N_2(\zeta)$. Hence, using the bounds on $N_1(\zeta)$ and $N_2(\zeta)$, (65) and the bound on the supremum of a Gaussian process in terms of the Dudley integral (see, e.g., Lifshits (1995), Section 14), we conclude that

$$\int_0^\sigma [\sqrt{\ln N_1(\zeta)} + \sqrt{\ln N_2(\zeta)}] \, d\zeta \leq c\sigma \sqrt{\ln \sigma}.$$

The first statement of the lemma is proved.

$2^0$. For the family of kernels $\mathcal{K}_{AH}$ of Section 4.3, we have $\nu = h \in \mathcal{H}_\gamma$, where $\mathcal{H}_\gamma$ is defined in (40). Note that the set $\Sigma_{\Theta_{AH}}$ consists of the single point

$$\sigma^* = \|G\|_1 \|G\|_2 [\varepsilon \sqrt{\ln \ln(1/\epsilon)}]^{-1/(2\gamma+1)}.$$

It follows from Lemma 2 that the semi-metric $\rho_{AH}$ of this process admits the following upper bound:

$$\rho_{AH}[(\mu, \nu), (\mu, \nu')] \leq 2\sigma^* M \left( \sum_{i=1}^d \left| 1 - \frac{h_i}{h'_i} \right|^2 \right)^{1/2}.$$

The number of balls $N(\zeta)$ of radius $\zeta$ in the above semi-metric covering the index set $\mathcal{H}_\gamma$ does not exceed

$$N(\zeta) \leq [\ln(h_{\max}/h_{\min}) \ln^{-1}(1 + c_1 \zeta/\sigma^*)]^d.$$

Hence, applying Lemma A.4 (see Appendix), we obtain

$$\int_0^{\sigma^*} \sqrt{\ln N(\zeta)} \, d\zeta \leq c_2 \sigma^* [\sqrt{\ln \ln(h_{\max}/h_{\min})} + 1].$$

$3^0$. For family of kernels $\mathcal{K}_B$, we have $\nu = h_1 \in [h_{\min}, h_{\max}]$ and

$$\tilde{\sigma}_{h_1} = \|G^*\|_1 \|G^*\|_2 h_1^{-d/2}, \qquad \sigma_{\min} = \|G^*\|_1 \|G^*\|_2 h_{\max}^{-d/2}, \qquad \sigma_{\max} = \|G^*\|_1 \|G^*\|_2 h_{\min}^{-d/2}.$$

According to Lemma 2,

$$\rho_B[(\mu, \nu), (\mu, \nu')] \leq 2\sigma \left\{ Md \left| 1 - \frac{h_1}{h'_1} \right| + 2 \left| 1 - \left( \frac{h'_1}{h_1} \right)^d \right| \right\}$$

for all $\nu = h_1$, $\nu' = h'_1$ such that $\tilde{\sigma}_{h_1} \vee \tilde{\sigma}_{h'_1} \leq \sigma$.

For fixed $\sigma \in [\sigma_{\min}, \sigma_{\max}]$, we set $h_\sigma = (\|G^*\|_1 \|G^*\|_2 \sigma^{-1})^{d/2}$. The number of balls $N(\zeta)$ of radius $\zeta$ in the semi-metric $\rho_B$ covering the set $[h_\sigma, h_{\max}]$ does not exceed

$$N(\zeta) \leq \ln(c_1 h_{\max} \sigma^{d/2}) \ln^{-1}(1 + c_2 [\zeta \sigma^{-1}]^{1/d}).$$

Hence,

$$\int_0^\sigma \sqrt{\ln N(\zeta)} \, d\zeta \leq c\sigma \sqrt{\ln(1 + \ln(\sigma/\sigma_{\min}))}. \qquad \square$$

### 6.2. Proof of Theorem 2

Throughout the proof, $c, c_1, c_2, \ldots$ stand for constants depending only on $d$, $G$ and $r$.

We show that $\Theta_F(\mathcal{K}_{SI})$ is non-empty for any $F \in \mathbb{F}_{SI}(\alpha, L)$. Assume that $F(t) = f(\omega_0^T t)$, $\omega_0 \in \mathbb{S}^{d-1}$, where $f \in \mathbb{H}_1(\alpha, L)$.

First, we note that in view of (36) and (38), there exist constants $c_1, c_2$ such that

$$c_1 \sqrt{\frac{1}{h_1} \ln \frac{1}{h_1}} \leq Q(\tilde{\sigma}_{h,E}) \leq c_2 \sqrt{\frac{1}{h_1} \ln \frac{1}{h_1}} \qquad \forall (h,E) \in \mathcal{H}_1 \times \mathcal{E}. \tag{66}$$

Consider the family of kernels $\mathcal{K}_{SI}^0 := \{G_{h,E} : h \in \mathcal{H}_1, E = E_0\} \subset \mathcal{K}_{SI}$, where $E_0$ is a fixed orthogonal matrix whose first column is $\omega_0$. Clearly, for any estimator associated with kernel $G_{h,E_0}$ from $\mathcal{K}_{SI}^0$, we have the following bound on the bias: $|B_{h,E_0}(x)| \leq L h_1^\alpha$ for all $x$. Moreover, by (17) and the fact that $M(\mathcal{K}_{SI}) = \|G\|_1$, we obtain

$$\tilde{B}_{h,E_0}(x) \leq \|G\|_1 \sup_y |B_{h,E_0}(y)| \leq \|G\|_1 L h_1^\alpha \qquad \forall h \in \mathcal{H}_1.$$

Let $h^* = (h_1^*, h_{\max}, \ldots, h_{\max})$ be defined by the balance equation

$$\|G\|_1 L (h_1^*)^\alpha = \tfrac{1}{2} \varepsilon Q(\tilde{\sigma}_{h^*, E_0}).$$

It then follows from (66) that

$$h_1^* = c_3 [(\varepsilon/L) \sqrt{\ln(\varepsilon/L)}]^{2/(2\alpha+1)}. \tag{67}$$

Note that for $\varepsilon$ small enough, $h_1^* \in [h_{\min}, h_{\max}]$ and by definition of $h_1^*$, $\tilde{B}_{h^*, E_0} \leq \tfrac{1}{2} \varepsilon Q(\tilde{\sigma}_{h^*, E_0})$.

We now show that $(h^*, E_0) \in \Theta_F(\mathcal{K}_{SI})$. To that end, fix $\sigma \in [\tilde{\sigma}_{h^*, E_0}, \sigma_{\max}]$. Consider the estimator associated with parameter $(h', E_0)$ such that $\tilde{\sigma}_{h', E_0} = \sigma$. Hence, by (36), $h_1' = c\sigma^{-2} \leq c\tilde{\sigma}_{h^*, E_0}^{-2} = h_1^*$ so that, in view of the monotonicity of the function $Q(\cdot)$,

$$\tilde{B}_{h', E_0} \leq \|G\|_1 h(h_1')^\alpha \leq \|G\|_1 L(h_1^*)^\alpha = \tfrac{1}{2} \varepsilon Q(\tilde{\sigma}_{h^*, E_0}) \leq \tfrac{1}{2} \varepsilon Q(\tilde{\sigma}_{h', E_0}).$$

This shows that $(h^*, E_0) \in \Theta_F(\mathcal{K}_{SI})$. Then, applying Theorem 1, we obtain

$$\{\mathbb{E}_F |\hat{F}_{SI}(x) - F(x)|^r\}^{1/r} \leq c \varepsilon Q(\tilde{\sigma}_{h^*, E_0}).$$

Substitution of (67) completes the proof.

### 6.3. Proof of Theorem 3

If $F \in \mathbb{F}_{AH}(\gamma, L)$, then there exists $\boldsymbol{\alpha}^* = (\alpha_1^*, \ldots, \alpha_d^*) \in \mathcal{A}_\gamma$ such that $F \in \mathbb{H}_d(\boldsymbol{\alpha}^*, L)$. Note that under the premise of the theorem, we have, for any $h \in \mathcal{H}_\gamma$, that

$$\varepsilon Q(\tilde{\sigma}_h) = c_1 \varepsilon \left( \prod_{i=1}^d h_i^{-1/2} \right) \sqrt{\ln \ln(h_{\max}/h_{\min})} = c_2 [\varepsilon \sqrt{\ln \ln(1/\varepsilon)}]^{2\gamma/(2\gamma+1)} = c_2 \varphi_\varepsilon.$$

Define $h^* = (h_1^*, \ldots, h_d^*)$ by the following relation:

$$d\|G\|_1 L(h_i^*)^{\alpha_i^*} = \frac{1}{2}\varepsilon Q(\tilde{\sigma}_{h^*}) = \frac{c_2}{2}\varphi_\varepsilon \qquad \forall i = 1, \ldots, d. \tag{68}$$

For $\varepsilon$ small enough, $h^* \in [h_{\min}, h_{\max}]^d$ and, clearly, $h^* \in \mathcal{H}_\gamma$. Let $\hat{F}_{h_*}$ be the estimator from $\mathcal{F}(\mathcal{K}_{AH})$ associated with kernel $G_{h^*}$ (see (33)). We have the following upper bound on the bias of this estimator: $\sup_x |B_{h^*}(x)| \leq L \sum_{i=1}^d (h_i^*)^{\alpha_i^*}$. Moreover, by (17),

$$\tilde{B}_{h^*}(x) \leq \|G_0\|_1 \sup_x |B_{h^*}(x)| \leq \|G\|_1 L \sum_{i=1}^d (h_i^*)^{\alpha_i^*} \leq d\|G\|_1 L(h_1^*)^{\alpha_1^*} = \tfrac{1}{2}\varepsilon Q(\tilde{\sigma}_{h^*}), \tag{69}$$

where the last inequality on the right-hand side follows from (68). Because the set $\Sigma_\Theta$ is the singleton $\{\tilde{\sigma}_{h^*}\}$, inequality (69) implies that $\Theta_F(\mathcal{K}_{AH})$ is non-empty. Application of Theorem 1 yields

$$\{\mathbb{E}_F|\hat{F}_{AH}(x) - F(x)|^r\}^{1/r} \leq c_3 \varepsilon Q(\tilde{\sigma}_{h^*}) = c_4 \varphi_\varepsilon.$$

The theorem is thus proved.

### 6.4. Proof of Theorem 4

Before turning to the proof of the theorem, let us make some remarks which will be used in the subsequent proof.

1. Let $s[q] = [s - d/p + d/q] \wedge s$, $q \in [1, \infty]$. Then, due to the inclusion theorem for Besov balls (Nikolskii (1975)), we have

$$\mathbb{B}_{p,\infty}^s(d, L) \subseteq \mathbb{B}_{q,\infty}^{s[q]}(d, L). \tag{70}$$

In particular,

$$\mathbb{B}_{p,\infty}^s(d, L) \subseteq \mathbb{B}_{\infty,\infty}^{s-d/p}(d, L) \subset \mathbb{C}(\mathcal{D}_0). \tag{71}$$

The last inclusion follows from the assumption of the theorem that $s - d/p > 0$. It also implies $s[q] > 0$.

2. Let us introduce the following notation. For any $\mathfrak{h} \in (0, h_{\max}]$, let $h(\mathfrak{h}) = (\mathfrak{h}, \ldots, \mathfrak{h}) \in \mathcal{H}_B$ and define

$$B_\mathfrak{h}(\cdot) := \int G_{h(\mathfrak{h})}^*(t - \cdot) F(t)\,\mathrm{d}t - F(\cdot) =: B_\mu(x), \qquad \mu = h(\mathfrak{h}),$$

$$B_{\mathfrak{h},\mathfrak{h}'}(\cdot) := B_{\mu,\nu}(\cdot) - B_\nu(x) = \int G_{h(\mathfrak{h}')}^*(y - \cdot) B_\mathfrak{h}(y)\,\mathrm{d}y, \qquad \nu = h(\mathfrak{h}'), \tag{72}$$

$$\tilde{B}_\mathfrak{h}(\cdot) := \sup_{\mathfrak{h}' \in (0, h_{\max}]} |B_{\mathfrak{h},\mathfrak{h}'}(\cdot) - B_{\mathfrak{h}'}(\cdot)| \vee |B_\mathfrak{h}(\cdot)|.$$

Let $\mathbb{K}_b(x), x \in \mathcal{D}_0$, be the cube centered at $x$ with side length equal to $b$. The possible values of $b$ are found from the condition $x + b \in \mathcal{D}$ for all $x \in \mathcal{D}_0$ and, later, $\sup_b$ denotes the supremum over this set.

Assuming, without loss of generality, that the support of the function $g$ belongs to $[-1/2(\lfloor s \rfloor + 2), 1/2(\lfloor s \rfloor + 2)]^d$ (it also implies that the support of $G^*$ belongs to $[-1/2, 1/2]^d$), we obtain from (72) that for all $x \in \mathcal{D}_0$ and $\mathfrak{h} \in (0, h_{\max}]$,

$$\tilde{B}_\mathfrak{h}(x) \leq C_1(s,g) \sup_b \frac{1}{b^d} \int_{\mathbb{K}_b(x)} |B_\mathfrak{h}(y)| \, \mathrm{d}y =: C_1(s,g) B_\mathfrak{h}^{(\max)}(x), \tag{73}$$

where $C_1(s,g)$ is a constant depending only on $\|g\|_\infty$ and $s$. To obtain (73), we used the fact that $\tilde{B}_\mathfrak{h}(x)$ is a continuous (even uniformly continuous) function of $\mathfrak{h}$, since $F$ is uniformly continuous on $\mathcal{D}_0$ and $G^*$ is bounded. The uniform continuity of $F$ follows from (71) and the compactness of $\mathcal{D}_0$. Note that $B_\mathfrak{h}^{(\max)}(\cdot)$ is the *Hardy–Littlewood maximal function* of $B_\mathfrak{h}(\cdot)$ (see, e.g., Wheeden and Zygmund (1977), Chapter 9, Section 3).

3. The operator $\Delta_a^l$ has the following representation:

$$\Delta_a^l F(x) = \sum_{j=0}^l \binom{l}{j} (-1)^{j+l} F(x + ja) \qquad \forall l \geq 1, \ \forall a > 0.$$

Therefore,

$$(-1)^{l+1} \Delta_a^l F(x) = \left[ \sum_{j=1}^l \binom{l}{j} (-1)^{j+1} F(x + ja) \right] - F(x).$$

Using this formula and the definition of the function $G^*$, we obtain, for any $\mathfrak{h} \in (0, h_{\max}]$,

$$\begin{aligned}
B_\mathfrak{h}(x) &= \int G^*(u)\{F(x + u\mathfrak{h}) - F(x)\} \, \mathrm{d}u \\
&= \sum_{j=1}^{\lfloor s \rfloor + 2} (-1)^{j+1} \binom{\lfloor s \rfloor + 2}{j} \frac{1}{j^d} \int g(u/j)\{F(x + u\mathfrak{h}) - F(x)\} \, \mathrm{d}u \\
&= \int g(v) \left\{ \sum_{j=1}^{\lfloor s \rfloor + 2} (-1)^{j+1} \binom{\lfloor s \rfloor + 2}{j} F(x + vj\mathfrak{h}) - F(x) \right\} \mathrm{d}v \\
&= (-1)^{\lfloor s \rfloor + 3} \int g(v) \Delta_{v\mathfrak{h}}^{\lfloor s \rfloor + 2} F(x) \, \mathrm{d}v.
\end{aligned}$$

Therefore,

$$|B_\mathfrak{h}(x)| \leq C_2(g) \int_{[-1/2, 1/2]^d} |\Delta_{v\mathfrak{h}}^{\lfloor s \rfloor + 2} F(x)| \, \mathrm{d}v =: C_2(g) \mathcal{B}_\mathfrak{h}(x) \qquad \forall x \in \mathcal{D}_0, \tag{74}$$

where $C_2(g)$ is a constant depending only on $\|g\|_\infty$. Here, we used the fact that the support of $g$ belongs to $[-1/2, 1/2]^d$.

Finally, from (73) and (74), we get

$$\tilde{B}_{\mathfrak{h}}(x) \leq C_3(s,g)\mathcal{B}_{\mathfrak{h}}^{(\max)}(x) \qquad \forall \mathfrak{h} \in (0, h_{\max}], \ \forall x \in \mathcal{D}_0, \tag{75}$$

where, as before, $\mathcal{B}_{\mathfrak{h}}^{(\max)}(\cdot)$ is the *Hardy–Littlewood maximal function* of $\mathcal{B}_{\mathfrak{h}}(\cdot)$ and $C_3(s,g) = C_1(s,g)C_2(g)$.

The next important property of $\mathcal{B}_{\mathfrak{h}}^{(\max)}(\cdot)$ follows from the definition:

$$\sup_{\mathfrak{h} \in [\tau/2,\tau]} \mathcal{B}_{\mathfrak{h}}^{(\max)}(x) \leq 2^d \mathcal{B}_{\tau}^{(\max)}(x) \qquad \forall \tau \in (0, h_{\max}], \ \forall x \in \mathcal{D}_0. \tag{76}$$

Indeed, for any $\mathfrak{h} \in [\tau/2, \tau]$, we have

$$\mathcal{B}_{\mathfrak{h}}(\cdot) := \mathfrak{h}^{-d} \int_{[-\mathfrak{h}/2, \mathfrak{h}/2]^d} |\Delta_u^{\lfloor s \rfloor + 2} F(\cdot)| \, \mathrm{d}u \leq 2^d \tau^{-d} \int_{[-\tau/2, \tau/2]^d} |\Delta_u^{\lfloor s \rfloor + 2} F(\cdot)| \, \mathrm{d}u$$
$$= 2^d \mathcal{B}_{\tau}(\cdot).$$

4. Since $\|G^*_{h(\mathfrak{h})}\|_2 = \|G^*\|_2 \mathfrak{h}^{-d/2} =: \sigma_{\mathfrak{h}}$, the majorant can be rewritten in the form

$$Q(\mathfrak{h}) := Q(\sigma_{\mathfrak{h}}) = \mathcal{C} h_{\max}^{-d/2} Q_d^*(h_{\max}/\mathfrak{h}),$$

where $Q_d^*(z) = z^{d/2}\sqrt{1 + d\ln z}, \ z \geq 1$.

Thus, for our particular problem, the set $\Theta_F(\mathcal{K}_B)$ is

$$\Theta_F(\mathcal{K}_B) = \Theta_F^x(\mathcal{K}_B) = \{\mathfrak{h} \in (0, h_{\max}] : \tilde{B}_{\mathfrak{h}'}(x) \leq \tfrac{1}{2}\varepsilon Q(\mathfrak{h}'), \forall \mathfrak{h}' \leq \mathfrak{h}\}, \qquad x \in \mathcal{D}_0.$$

Note that (75) and (71) imply that for all $F \in \mathbb{B}_{p,\infty}^s(d, L)$ and all $\mathfrak{h} \in (0, h_{\max}]$,

$$\|\tilde{B}_{\mathfrak{h}}\|_\infty \leq C_3(s,g)\|\mathcal{B}_{\mathfrak{h}}^{(\max)}\|_\infty = C_3(s,g)\|\mathcal{B}_{\mathfrak{h}}\|_\infty \leq 2^d C_3(s,g) \sup_{v \in [-1,1]^d} \|\Delta_{v\mathfrak{h}}^{\lfloor s \rfloor + 2} F\|_\infty$$

$$\leq 2^d C_3(s,g)[\mathfrak{h}\sqrt{d}]^{s-d/p} \sup_a |a|_2^{-s+d/p}\|\Delta_a^{\lfloor s \rfloor + 2} F\|_\infty \leq LC_4(s,g,p,d)\mathfrak{h}^{s-d/p}.$$

Therefore, there exists a constant $c$, depending only on $s, g, p, d$ and $L$, such that

$$(0, \varepsilon^c] \subset \Theta_F^x(\mathcal{K}_B) \qquad \forall F \in \mathbb{B}_{p,\infty}^s(d, L), \ \forall x \in \mathcal{D}_0. \tag{77}$$

Putting, for all $F \in \mathbb{B}_{p,\infty}^s(d, L)$ and all $x \in \mathcal{D}_0$,

$$\mathfrak{h}_F(x) = \sup\{\mathfrak{h} : \mathfrak{h} \in \Theta_F^x(\mathcal{K}_B)\},$$

we obtain from Theorem 1 that

$$\mathcal{R}_{\mathbb{L}_r}(\tilde{F}) \leq C^r \varepsilon^r \sup_{F \in \mathbb{B}_{p,\infty}^s(d,L)} \int_{\mathcal{D}_0} Q^r(\mathfrak{h}_F(x)) \, \mathrm{d}x.$$

**Proof of Theorem 4.** As already mentioned, the function $\tilde{B}_{\mathfrak{h}}(\cdot)$ is continuous in $\mathfrak{h}$. Evidently, the function $Q(\mathfrak{h})$ is also continuous. Therefore, in view of the definition of $\mathfrak{h}_F(x)$, we have

$$\tilde{B}_{\mathfrak{h}_F(x)}(x) = \tfrac{1}{2}\varepsilon Q(\mathfrak{h}_F(x)) \qquad \forall x \in \mathcal{D}_0 : \mathfrak{h}_F(x) < h_{\max}. \tag{78}$$

Let $k_{\max} \in \mathbb{N}^*$ be chosen in such a way that $2^{-k_{\max}} h_{\max} \leq h_{\min} < 2^{1-k_{\max}} h_{\max}$, where $h_{\min} = \varepsilon^c$. Set

$$\Gamma_0 = \{x \in \mathcal{D}_0 : \mathfrak{h}_F(x) = h_{\max}\},$$
$$\Gamma_k = \{x \in \mathcal{D}_0 : 2^{-k} h_{\max} \leq \mathfrak{h}_F(x) < 2^{1-k} h_{\max}\}, \qquad k = \overline{1, k_{\max}}.$$

Note that the sets $(\Gamma_k, k = \overline{1, k_{\max}})$ form the partition of $\mathcal{D}_0$ since $\mathfrak{h}_F(x) \geq h_{\min}$ for all $x \in \mathcal{D}_0$, in view of (77). Therefore,

$$I(F) := \varepsilon^r \int_{\mathcal{D}_0} Q^r(\mathfrak{h}_F(x))\,\mathrm{d}x = \varepsilon^r \sum_{k=0}^{k_{\max}} \int_{\Gamma_k} Q^r(\mathfrak{h}_F(x))\,\mathrm{d}x =: \sum_{k=0}^{k_{\max}} I_k(F). \tag{79}$$

Let $q_k \in (1, r], k \in \mathbb{N}^*$, be a sequence of real numbers, to be specified later. Then, in view of (78), we get $\forall k = \overline{1, k_{\max}}$,

$$I_k(F) = \varepsilon^{r-q_k} 2^{q_k} \int_{\Gamma_k} Q^{r-q_k}(\mathfrak{h}_F(x))(\tilde{B}_{\mathfrak{h}_F(x)}(x))^{q_k}\,\mathrm{d}x. \tag{80}$$

It follows from (75) and (76) that

$$\tilde{B}_{\mathfrak{h}_F(x)}(x) \leq 2^d C_3(s,g) \mathcal{B}^{(\max)}_{2^{1-k}h_{\max}}(x) \qquad \forall x \in \Gamma_k. \tag{81}$$

Moreover,

$$Q(\mathfrak{h}_F(x)) \leq Q(2^{-k} h_{\max}) = \mathcal{C} h_{\max}^{-d/2} 2^{kd/2} \sqrt{1 + kd\ln 2} \qquad \forall x \in \Gamma_k. \tag{82}$$

Thus, we have, from (80), (81) and (82), that

$$I_k(F) \leq C_1 [\varepsilon h_{\max}^{-d/2}]^{r-q_k} 2^{kd(r-q_k)/2} k^{(r-q_k)/2} \|\mathcal{B}^{(\max)}_{2^{1-k}h_{\max}}\|_{q_k}^{q_k}. \tag{83}$$

Here and later, we denote by $C_1, C_2, \ldots$, the constants depending on $d, s, p, r, g, L$, but independent of $F$ and $\varepsilon$.

We have, for all $\mathfrak{h} \in (0, h_{\max}]$ and all $F \in \mathbb{B}_{p,\infty}^s(d, L)$,

$$\|\mathcal{B}_{\mathfrak{h}}^{(\max)}\|_{q_k}^{q_k} \leq C(q_k) \|\mathcal{B}_{\mathfrak{h}}\|_{q_k}^{q_k} = C(q_k) \left\| \int_{[-1/2, 1/2]^d} |\Delta_{v\mathfrak{h}}^{\lfloor s \rfloor + 2} F|\,\mathrm{d}v \right\|_{q_k}^{q_k}$$
$$\leq C(q_k) \left[ \int_{[-1/2, 1/2]^d} \|\Delta_{v\mathfrak{h}}^{\lfloor s \rfloor + 2} F\|_{q_k}\,\mathrm{d}v \right]^{q_k}$$

(84)
$$\leq C_2 \left[ \mathfrak{h}^{s[q_k]} \sup_a |a|_2^{-s[q_k]} \|\Delta_a^{\lfloor s \rfloor + 2} F\|_{q_k} \right]^{q_k}$$
$$\leq (LC_2 \mathfrak{h}^{s[q_k]})^{q_k} = C_3 \mathfrak{h}^{q_k s[q_k]}.$$

Let us comment on the proof of (84). The first inequality follows from (Wheeden and Zygmund (1977), Theorem 9.16), where the constant $C(q_k)$ depends only on $q_k$ and, moreover, $\sup_{1 \leq q \leq r} C(q) < \infty$ for any fixed $r$. The second inequality follows from the Minkowski inequality for integrals. The last inequality is a consequence of (70).

Substituting $\mathfrak{h} = 2^{1-k} h_{\max}$ in (84), we finally obtain, from (79) and (83), that for any $F \in \mathbb{B}_{p,\infty}^s(d, L)$,

$$I(F) \leq C_4 \left[ \left( \frac{\varepsilon}{h_{\max}^{d/2}} \right)^r + \sum_{k=1}^{k_{\max}} (\varepsilon h_{\max}^{-d/2})^{r-q_k} (h_{\max})^{q_k s[q_k]} 2^{-k\lambda_k} k^{(r-q_k)/2} \right],$$
$$\lambda_k = q_k s[q_k] - (r - q_k) d/2, \qquad k = \overline{1, k_{\max}}.$$
(85)

Let us now consider three cases.

*Case 1.* $sp > \frac{d(r-p)}{2}$. Choose $q_k = r \wedge p$ for all $k = \overline{1, k_{\max}}$ and recall that $h_{\max} = \varepsilon^{2/(s+d/2)}$. Therefore,

$$s[q_k] = s, \qquad \lambda_k = \lambda := s(r \wedge p) - \frac{d(r - r \wedge p)}{2} > 0.$$

Moreover,

$$(h_{\max})^s = \frac{\varepsilon}{h_{\max}^{d/2}} = \varphi_\varepsilon.$$

Thus, we obtain from (85), for any $F \in \mathbb{B}_{p,\infty}^s(d, L)$,

$$I(F) \leq C_4 \left[ \varphi_\varepsilon^r + \varphi_\varepsilon^r \sum_{k=1}^{\infty} 2^{-k\lambda} k^{(r-r\wedge p)/2} \right] \leq C_5 \varphi_\varepsilon^r.$$

*Case 2.* $sp = \frac{d(r-p)}{2}$. Choose $q_k = p$ for all $k = \overline{1, k_{\max}}$ and recall that $h_{\max} = 1/2$. Therefore,

$$s[q_k] = s, \qquad \lambda_k = 0.$$

Taking into account that $k_{\max} \sim \ln(1/\varepsilon)$, we obtain from (85) that for any $F \in \mathbb{B}_{p,\infty}^s(d, L)$,

$$I(F) \leq C_6[\varepsilon^r + \varepsilon^{r-p}[\ln(1/\varepsilon)]^{(r-p)/2+1}] \leq C_7 \varphi_\varepsilon^r.$$

The last inequality follows from the relation $r - p = \frac{2rs}{2s+d}$.

*Case 3.* $sp < \frac{d(r-p)}{2}$. Recall that $h_{\max} = 1/2$ and $\varphi_\varepsilon = [\varepsilon\sqrt{\ln 1/\varepsilon}]^{s-d(1/p-1/r)/(s-d(1/p-1/2))}$. Let $\mathfrak{h}_\varepsilon = [\varepsilon\sqrt{\ln 1/\varepsilon}]^z$, $z = 1/(s - d/p + d/2)$, and define

$$q_k = \begin{cases} p, & \text{if } 1 \leq k < k^*, \\ r, & \text{if } k \geq k^*, \end{cases}$$

where $k^* \in \mathbb{N}^*$ is chosen from the relation $2^{-(k^*+1)} < \mathfrak{h}_\varepsilon \leq 2^{-k^*}$. Noting that

$$\lambda_k = \begin{cases} \lambda_1 := sp - \frac{d(r-p)}{2} < 0, & \text{if } 1 \leq k \leq k^*, \\ \lambda_2 := (s - d/p + d/r)r, & \text{if } k \geq k^* + 1 \end{cases}$$

and again taking into account that $k_{\max} \sim \ln(1/\varepsilon)$, we get, from (85), that for any $F \in \mathbb{B}^s_{p,\infty}(d, L)$,

$$I(F) \leq C_8 \left[ \varepsilon^r + [\varepsilon\sqrt{\ln 1/\varepsilon}]^{r-p} \sum_{k=1}^{k^*} 2^{-k\lambda_1} + \sum_{k=k^*+1}^{\infty} 2^{-k\lambda_2} \right]$$

$$\leq C_9 [[\varepsilon\sqrt{\ln 1/\varepsilon}]^{r-p} 2^{-k^*\lambda_1} + 2^{-(k^*+1)\lambda_2}]$$

$$\leq C_9 [[\varepsilon\sqrt{\ln 1/\varepsilon}]^{r-p} \mathfrak{h}_\varepsilon^{\lambda_1} + \mathfrak{h}_\varepsilon^{\lambda_2}].$$

To obtain the second inequality, we used the fact that $\lambda_1 < 0$ and $\lambda_2 > 0$. It remains to note that, in view of the definition of $\mathfrak{h}_\varepsilon$,

$$[\varepsilon\sqrt{\ln 1/\varepsilon}]^{r-p} \mathfrak{h}_\varepsilon^{\lambda_1} = \mathfrak{h}_\varepsilon^{\lambda_2} = \varphi_\varepsilon^r. \qquad \square$$

## Appendix

**Proof of Proposition 1.** We have, for all $\mu, \nu \in \Theta$ and all $x \in \mathcal{D}_0$,

$$B_{\mu,\nu}(x) = \int_{\mathcal{D}} K_{\mu,\nu}(t,x) F(t) \, dt - F(x) = \int_{\mathcal{D}} \left[ \int_{\mathcal{D}_1} K_\mu(t,y) K_\nu(y,x) \, dy \right] F(t) \, dt - F(x)$$

$$= \int_{\mathcal{D}_1} K_\nu(y,x) \left[ \int_{\mathcal{D}} K_\mu(t,y) F(t) \, dt \right] dy - F(x)$$

$$= \int_{\mathcal{D}_1} K_\nu(y,x) \left[ \int_{\mathcal{D}} K_\mu(t,y) \{F(t) - F(y) + F(y)\} \, dt \right] dy - F(x)$$

$$= \int_{\mathcal{D}_1} K_\nu(y,x) B_\mu(y) \, dy + \int_{\mathcal{D}_1} K_\nu(y,x) F(y) \, dy - F(x)$$

$$= \int_{\mathcal{D}} K_\nu(y,x) B_\mu(y) \, dy + \int_{\mathcal{D}} K_\nu(y,x) F(y) \, dy - F(x)$$

$$= \int_{\mathcal{D}} K_\nu(y,x) B_\mu(y) \, dy + B_\nu(x).$$

The fifth and sixth equalities follow from the second and first lines of (7), respectively. □

**Lemma A.1.** *For any $\mu \in \Theta$ and $r > 0$, we have*

$$\left\{\mathbb{E} \sup_{\nu:\tilde{\sigma}_\nu \leq \tilde{\sigma}_\mu} |\xi_{\nu,\mu} - \xi_\nu|^r\right\}^{1/r} \leq C_r\{e_{\mathcal{K}_\Theta}(\tilde{\sigma}_\mu) + 2\tilde{\sigma}_\mu\},$$

*where*

$$C_r = \begin{cases} 1, & r \leq 1, \\ \left[8r \int_0^\infty (t \vee 1)^{r-1} \exp\left(-\frac{t^2}{2}\right) dt\right]^{1/r}, & r > 1. \end{cases}$$

**Proof.** For brevity in the proof below, we will write $e(\cdot) = e_{\mathcal{K}_\Theta}(\cdot)$.

The statement for $r \leq 1$ follows immediately from the definition of the function $e(\cdot)$. If $r > 1$, then

$$\mathbb{E} \sup_{\nu:\tilde{\sigma}_\nu \leq \tilde{\sigma}_\mu} |\xi_{\nu,\mu} - \xi_\nu|^r = r \int_0^\infty t^{r-1} \mathbb{P}\left\{\sup_{\nu:\tilde{\sigma}_\nu \leq \tilde{\sigma}_\mu} |\xi_{\nu,\mu} - \xi_\nu| > t\right\} dt$$

$$\leq e^r(\tilde{\sigma}_\mu) + r \int_{e(\tilde{\sigma}_\mu)}^\infty t^{r-1} \mathbb{P}\left\{\sup_{\nu:\tilde{\sigma}_\nu \leq \tilde{\sigma}_\mu} |\xi_{\nu,\mu} - \xi_\nu| > t\right\} dt$$

$$= e^r(\tilde{\sigma}_\mu) + r \int_0^\infty [t + e(\tilde{\sigma}_\mu)]^{r-1} \mathbb{P}\left\{\sup_{\nu:\tilde{\sigma}_\nu \leq \tilde{\sigma}_\mu} |\xi_{\nu,\mu} - \xi_\nu| - e(\tilde{\sigma}_\mu) > t\right\} dt$$

$$\leq e^r(\tilde{\sigma}_\mu) + 2r \int_0^\infty [t + e(\tilde{\sigma}_\mu)]^{r-1} \exp\left\{-\frac{t^2}{2\sup_{\nu:\tilde{\sigma}_\nu \leq \tilde{\sigma}_\mu} \sigma_{\mu,\nu}^2}\right\} dt,$$

where the last inequality follows from the fact that $e(\sigma) \geq e_0(\sigma)$ and Lemma A.3 below; recall that $\text{var}(\xi_{\nu,\mu} - \xi_\nu) = \sigma_{\mu,\nu}^2$. Inequality (18) implies that $\sup_{\nu:\tilde{\sigma}_\nu \leq \tilde{\sigma}_\mu} \sigma_{\mu,\nu} \leq 2\tilde{\sigma}_\mu$; hence, continuing the preceding chain of inequalities, we obtain

$$\leq e^r(\tilde{\sigma}_\mu) + 2r \int_0^\infty [t + e(\tilde{\sigma}_\mu)]^{r-1} \exp\left\{-\frac{t^2}{8\tilde{\sigma}_\mu^2}\right\} dt$$

$$= e^r(\tilde{\sigma}_\mu) + 4r\tilde{\sigma}_\mu \int_0^\infty \{2t\tilde{\sigma}_\mu + e(\tilde{\sigma}_\mu)\}^{r-1} \exp(-t^2/2) \, dt$$

$$\leq e^r(\tilde{\sigma}_\mu) + 4r\tilde{\sigma}_\mu [2\tilde{\sigma}_\mu + e(\tilde{\sigma}_\mu)]^{r-1} \int_0^\infty (t \vee 1)^{r-1} \exp(-t^2/2) \, dt$$

$$\leq 8r[2\tilde{\sigma}_\mu + e(\tilde{\sigma}_\mu)]^r \int_0^\infty (t \vee 1)^{r-1} \exp(-t^2/2) \, dt.$$

This completes the proof. □

**Lemma A.2.** *Let assumption* (E) *hold and let the function $Q$ be given by (21) with $\varkappa_0 \geq 2C_e$ and $\varkappa_1 \geq 64$. Then, for any $\mu \in \Theta$ and $t \geq 0$,*

$$\mathbb{P}\left\{\sup_{\nu:\tilde{\sigma}_\nu \geq \tilde{\sigma}_\mu} \left[|\xi_{\mu,\nu} - \xi_\nu| - \frac{1}{2}Q(\tilde{\sigma}_\nu)\right] > t\right\} \leq 4\left(\frac{\sigma_{\min}}{\tilde{\sigma}_\mu}\right)^{\varkappa_1/64} \exp\left\{-\frac{t^2}{16\tilde{\sigma}_\mu^2}\right\}. \tag{86}$$

*Moreover, if $\varkappa_1 \geq 128r$, then*

$$\left\{\mathbb{E}\left[\sup_{\nu:\tilde{\sigma}_\nu \geq \tilde{\sigma}_\mu} |\xi_{\mu,\nu} - \xi_\nu| - \tfrac{1}{2}Q(\tilde{\sigma}_\nu)\right]_+^r\right\}^{1/r} \leq C\sigma_{\min}, \tag{87}$$

*where $C$ is a constant depending only on $r$.*

**Proof.** As previously, we will write $e(\cdot) = e_{\mathcal{K}_\Theta}(\cdot)$.

Define $N_k = \{\nu : 2^{k-1}\tilde{\sigma}_\mu \leq \tilde{\sigma}_\nu < 2^k \tilde{\sigma}_\mu\}$ for $k = 1, 2, \ldots$ and write

$$\mathbb{P}\left\{\sup_{\nu:\tilde{\sigma}_\nu \geq \tilde{\sigma}_\mu} [|\xi_{\mu,\nu} - \xi_\nu| - \tfrac{1}{2}Q(\tilde{\sigma}_\nu)] > t\right\} \leq \sum_{k=1}^{\infty} \mathbb{P}\left\{\sup_{\nu \in N_k} [|\xi_{\mu,\nu} - \xi_\nu| - \tfrac{1}{2}Q(\tilde{\sigma}_\nu)] > t\right\}. \tag{88}$$

Since $Q(\sigma)$ is monotone increasing in $\sigma$,

$$\mathbb{P}\left\{\sup_{\nu \in N_k}\left[|\xi_{\mu,\nu} - \xi_\nu| - \frac{1}{2}Q(\tilde{\sigma}_\nu)\right] > t\right\}$$

$$\leq \mathbb{P}\left\{\sup_{\nu \in N_k} |\xi_{\mu,\nu} - \xi_\nu| > t + \frac{1}{2}Q(2^{k-1}\tilde{\sigma}_\mu)\right\}$$

$$\leq \mathbb{P}\left\{\sup_{\nu:\tilde{\sigma}_\nu \leq 2^k\tilde{\sigma}_\mu} |\xi_{\mu,\nu} - \xi_\nu| - e(2^k\tilde{\sigma}_\mu) > t + \frac{1}{2}Q(2^{k-1}\tilde{\sigma}_\mu) - e(2^k\tilde{\sigma}_\mu)\right\}$$

$$= \mathbb{P}\bigg\{\sup_{\nu:\tilde{\sigma}_\nu \leq 2^k\tilde{\sigma}_\mu} |\xi_{\mu,\nu} - \xi_\nu| - e(2^k\tilde{\sigma}_\mu) > t + \frac{1}{2}\varkappa_0 e(2^{k-1}\tilde{\sigma}_\mu)$$

$$+ 2^{k-2}\tilde{\sigma}_\mu\sqrt{1 + \varkappa_1 \ln \frac{2^{k-1}\tilde{\sigma}_\mu}{\sigma_{\min}}} - e(2^k\tilde{\sigma}_\mu)\bigg\}$$

$$\leq \mathbb{P}\left\{\sup_{\nu:\tilde{\sigma}_\nu \leq 2^k\tilde{\sigma}_\mu} |\xi_{\mu,\nu} - \xi_\nu| - e(2^k\tilde{\sigma}_\mu) > t + 2^{k-2}\tilde{\sigma}_\mu\sqrt{1 + \varkappa_1 \ln \frac{2^{k-1}\tilde{\sigma}_\mu}{\sigma_{\min}}}\right\},$$

where the last inequality follows by assumption (E) and choice of the constant $\varkappa_0$. By (18),

$$\sup_{\nu:\tilde{\sigma}_\nu \leq 2^k\tilde{\sigma}_\mu} \mathrm{var}(\xi_{\mu,\nu} - \xi_\nu) = \sup_{\nu:\tilde{\sigma}_\nu \leq 2^k\tilde{\sigma}_\mu} \sigma_{\mu,\nu}^2 \leq 2^{2k+1}\tilde{\sigma}_\mu^2;$$

hence, using Lemma A.3, we obtain

$$\mathbb{P}\left\{\sup_{\nu \in N_k}\left[|\xi_{\mu,\nu} - \xi_\nu| - \frac{1}{2}Q(\sigma_\nu)\right] > t\right\} \le 2\exp\left\{-\frac{1}{2}\frac{(t+a_k)^2}{b_k^2}\right\},$$

where we have denoted for brevity

$$a_k = 2^{k-2}\tilde{\sigma}_\mu\sqrt{1 + \varkappa_1\ln(2^{k-1}\tilde{\sigma}_\mu/\sigma_{\min})}, \qquad b_k = 2^{k+1/2}\tilde{\sigma}_\mu.$$

Noting that

$$\exp\left\{-\frac{(t+a_k)^2}{2b_k^2}\right\} \le \exp\left\{-\frac{t^2}{2b_k^2}\right\}\exp\left\{-\frac{a_k^2}{2b_k^2}\right\}$$

$$\le \exp\left\{-\frac{t^2}{16\tilde{\sigma}_\mu^2}\right\}2^{-(k-1)\varkappa_1/64}\left(\frac{\sigma_{\min}}{\tilde{\sigma}_\mu}\right)^{\varkappa_1/64},$$

we have

$$\mathbb{P}\left\{\sup_{\nu \in N_k}\left[|\xi_{\mu,\nu} - \xi_\nu| - \frac{1}{2}Q(\sigma_\nu)\right] > t\right\} \le 2^{1-(k-1)\varkappa_1/64}\exp\left\{-\frac{t^2}{16\tilde{\sigma}_\mu^2}\right\}\left(\frac{\sigma_{\min}}{\tilde{\sigma}_\mu}\right)^{\varkappa_1/64}. \quad (89)$$

Summing up over $k = 1, 2, \ldots$ and taking into account (88), we arrive at (86).

We now prove (87). Using (89), we have

$$\mathbb{E}\left[\sup_{\nu \in N_k}|\xi_{\mu,\nu} - \xi_\nu| - \frac{1}{2}Q(\tilde{\sigma}_\nu)\right]_+^r \le 2^{1-(k-1)\varkappa_1/64}\left(\frac{\sigma_{\min}}{\tilde{\sigma}_\mu}\right)^{\varkappa_1/64} r\int_0^\infty t^{r-1}\exp\left\{-\frac{t^2}{16\tilde{\sigma}_\mu^2}\right\}\mathrm{d}t$$

$$\le c2^{-(k-1)\varkappa_1/64}\left(\frac{\sigma_{\min}}{\tilde{\sigma}_\mu}\right)^{\varkappa_1/64}\tilde{\sigma}_\mu^r.$$

This implies that

$$\left\{\mathbb{E}\left[\sup_{\nu:\sigma_\nu \ge \sigma_\mu}|\xi_{\mu,\nu} - \xi_\nu| - \frac{1}{2}Q(\tilde{\sigma}_\nu)\right]_+^r\right\}^{1/r} \le \left\{\sum_{k=1}^\infty \mathbb{E}\left[\sup_{\nu \in N_k}|\xi_{\mu,\nu} - \xi_\nu| - \frac{1}{2}Q(\tilde{\sigma}_\nu)\right]_+^r\right\}^{1/r}$$

$$\le c\tilde{\sigma}_\mu\left(\frac{\sigma_{\min}}{\tilde{\sigma}_\mu}\right)^{\varkappa_1/(64r)}\left\{\sum_{k=0}^\infty 2^{-k\varkappa_1/64}\right\}^{1/r}$$

$$\le c\tilde{\sigma}_\mu\left(\frac{\sigma_{\min}}{\tilde{\sigma}_\mu}\right)^2 \le c\sigma_{\min},$$

because of our choice of $\varkappa_1$. $\square$

The next result can be found in, for example, Adler and Taylor (2007).

**Lemma A.3 (Borell, Tsirelson and Sudakov).** *Let $X_t$, $t \in T$, be a centered Gaussian process, a.s. bounded on $T$. Then, for all $u > 0$,*

$$\mathbb{P}\left\{\sup_{t \in T} X_t > \mathbb{E} \sup_{t \in T} X_t + u\right\} \leq \exp\{-u^2/2\sigma_T^2\}$$

*and hence*

$$\mathbb{P}\left\{\sup_{t \in T} |X_t| > \mathbb{E} \sup_{t \in T} |X_t| + u\right\} \leq 2\exp\{-u^2/2\sigma_T^2\},$$

*where $\sigma_T^2 = \sup_{t \in T} \mathrm{var}(X_t)$.*

**Lemma A.4.** *Let $a > 0$ and $a\sigma < \exp(1) - 1$. Then,*

$$\int_0^\sigma \sqrt{\ln \ln^{-1}(1 + \eta a)}\, \mathrm{d}\eta \leq \frac{\exp(1)}{a} \frac{\sqrt{\ln \ln^{-1}(1 + a\sigma)}}{\ln(1 + a\sigma)} \left\{1 + \frac{1}{2 \ln \ln^{-1}(1 + a\sigma)}\right\}.$$

The proof is immediate.

## Acknowledgment

The first author was supported by an ISF research grant.